\documentclass{birkmult}

\usepackage[all]{xy}

\theoremstyle{plain}
  \newtheorem{theorem}{Theorem}
  
  \newtheorem{conjecture}{Conjecture}
  
  \newtheorem{proposition}{Proposition}
  \newtheorem{observation}{Observation}
\theoremstyle{definition}
  \newtheorem{definition}{Definition}
\theoremstyle{remark}
  \newtheorem{remark}{Remark}

\newcommand{\an}{\text{\rm an}}

\newcommand\even{\text{\rm even}}
\newcommand\odd{\text{\rm odd}}

\newcommand\ac{\text{\rm ac}}

\newcommand\bb{\mathfrak b}

\newcommand\id{\text{\rm id}}
\newcommand{\Z}{\mathbb Z}

\newcommand{\R}{\mathbb R}
\newcommand{\C}{\mathbb C}

\newcommand{\Eul}{\mathfrak{Eul}}
\newcommand\vf{\mathfrak X}
\newcommand{\e}{\mathfrak e}
\newcommand\Or{\mathcal O}

\newcommand\comb{\text{\rm comb}}

\DeclareMathOperator{\llet}{let}
\DeclareMathOperator{\GL}{GL}

\DeclareMathOperator\eend{end}

\DeclareMathOperator\Tor{Tor}
\DeclareMathOperator\img{img}
\DeclareMathOperator\Emb{Emb}
\DeclareMathOperator\Gr{Gr}
\DeclareMathOperator\tr{tr}
\DeclareMathOperator\inter{int}
\DeclareMathOperator\pt{pt}
\DeclareMathOperator{\ind}{ind}

\DeclareMathOperator{\sgn}{sgn}

\DeclareMathOperator{\IND}{IND}
\DeclareMathOperator{\Vol}{vol}

\DeclareMathOperator{\Spect}{Spect}

\DeclareMathOperator{\Int}{Int}
\DeclareMathOperator{\grad}{grad}

\DeclareMathOperator{\cs}{cs}
\DeclareMathOperator{\ec}{e}
\DeclareMathOperator\PD{PD}
\DeclareMathOperator\Rep{Rep}

\renewcommand{\theenumi}{\roman{enumi}}
\renewcommand{\labelenumi}{(\theenumi)}

\begin{document}

\title[Torsions,   as a functions on the  space of representations]
      {Torsion,   as a function on the  space of\\ representations}

\author{Dan Burghelea}

\address{%Dan Burghelea,
         Dept. of Mathematics, The Ohio State University,\\
         231 West 18th Avenue, Columbus, OH 43210, USA.}

\email{burghele@mps.ohio-state.edu}

\author{Stefan Haller}

\address{%Stefan Haller,
         Department of Mathematics, University of Vienna,\\
         Nordbergstrasse 15, A-1090, Vienna, Austria.}

\email{stefan.haller@univie.ac.at}

\thanks{Part of this work was done while both authors enjoyed the
        hospitality of the Max Planck Institute for Mathematics in Bonn.
        A previous version was written while the second author enjoyed the
        hospitality of the Ohio State University.
        The second author was partially supported by the \emph{Fonds zur F\"orderung der
        wissenschaftlichen Forschung} (Austrian Science Fund),
        project number {\tt P17108-N04}}


%\subjclass[2000]{57R20, 58J52}

%\date{\today}

%%%%%%%%%%%%%%%%%%%%%%%%%%%%%%%%%%%%%%%%%%%%%%%%%%%%%%%%%%%%

\begin{abstract}
Riemannian Geometry, Topology and 
Dynamics permit to introduce partially defined holomorphic 
functions on the variety of representations of the fundamental group of a manifold.
The functions we consider are the complex valued Ray--Singer torsion, the Milnor--Turaev 
torsion, and the dynamical torsion. They are associated essentially to a 
closed smooth manifold equipped with a (co)Euler structure and a Riemannian metric 
in the first case, a smooth triangulation in the second case, and a smooth 
flow of type described in section~\ref{S:2} in the third case.
In this paper we define these functions, describe 
some of their properties and calculate them in some case. 
We conjecture that they are essentially equal and  
have analytic continuation to rational functions on the variety of representations. 
We discuss the case of one dimensional representations 
and other relevant situations when the conjecture is true.  
As particular cases of our torsions, we recognize familiar rational functions in 
topology such as the Lefschetz zeta function of a diffeomorphism, the dynamical 
zeta function of closed trajectories, and the Alexander polynomial of a knot. 
A numerical invariant derived from Ray--Singer torsion and associated to 
two homotopic acyclic representations is discussed in the last section.
\end{abstract}

\maketitle

\newpage

\setcounter{tocdepth}{1}
\tableofcontents

\section{Introduction}\label{S:intro}

For a finitely presented group $\Gamma$ denote by $\Rep(\Gamma;V)$ 
the algebraic set of all complex representations of $\Gamma$ on
the complex vector space $V$.
For a closed base pointed manifold $(M,x_0)$ with $\Gamma=\pi_1(M,x_0)$ 
denote by $\Rep^M(\Gamma;V)$ the algebraic closure of $\Rep^M_0(\Gamma;V)$,
the Zariski open set of representations $\rho\in\Rep(\Gamma;V)$ so that 
$H^*(M;\rho)=0$.
The manifold $M$ is called $V$-acyclic iff $\Rep^M(\Gamma;V)$, or
equivalently $\Rep^M_0(\Gamma;V)$, is non-empty. If $M$ is $V$-acyclic
then the Euler--Poincar\'e characteristic $\chi(M)$ vanishes. There are plenty of
$V$-acyclic manifolds.

If $\dim V=1$ then $\Rep(\Gamma;V)=(\mathbb C\setminus 0)^k\times F$, where $k$ denotes
the first Betty number of $M$, and $F$ is a finite Abelian group. If in addition $M$ 
is $V$-acyclic and $H_1(M;\mathbb Z)$ is torsion free, then 
$\Rep^M(\Gamma;V)=(\mathbb C\setminus 0)^k$. 
There are plenty of $V$-acyclic ($\dim V=1$) manifolds $M$ with $H_1(M;\mathbb Z)$ 
torsion free.

In this paper, to a $V$-acyclic manifold and an Euler or coEuler structure
we associate three partially defined holomorphic functions on
$\Rep^M(\Gamma;V)$, the complex valued Ray--Singer torsion, the
Milnor--Turaev torsion, and the dynamical torsion, and describe some of
their properties.  They are defined with the help of a Riemannian metric,
resp.\ smooth triangulation resp.\ a vector field with the properties listed
in section~\ref{SS:2.3}, but are independent of these data.

We conjecture that they are essentially equal and have analytic continuation to rational
functions on $\Rep^M(\Gamma;V)$ and discuss the cases when we know that this is true.
If $\dim V=1$ they are genuine rational functions of $k$ variables.

We calculate them in some cases and recognize familiar rational functions in topology
(Lefschetz zeta function of a diffeomorphism, dynamical zeta function of some flows, 
Alexander polynomial of a knot) as particular cases of our torsions, cf.\ section~\ref{S:7}.

The results answer the question
\begin{itemize}
\item[(Q)]
\emph{Is the Ray--Singer torsion the absolute value of a  holomorphic 
function on the space of representations?}\footnote{A similar question was 
considered in \cite{Q85} and a positive answer provided.}
\end{itemize} 
(for a related result consult \cite{BK05})
and  establish the analytic continuation of the dynamical torsion.
Both issues are  subtle when formulated inside the field of spectral geometry  
or of dynamical systems and can hardly be decided using internal technologies 
in these fields. There are  interesting dynamical implications on the growth 
of the number of instantons  and of closed trajectories, some of them improving 
on a conjecture formulated by S.P.~Novikov about the gradients of closed Morse one 
form, cf.\ section~\ref{S:8}.

\vskip .1in

This paper surveys results from \cite{BH04}, \cite{BH03'}, \cite{BH05} and reports 
on additional work in progress on these lines. Its contents is the following.

In section~\ref{S:2}, for the reader's convenience, we recall some less 
familiar characteristic forms used in this paper and describe the class of vector 
fields we use to define the dynamical torsion. These vector fields  have finitely 
many rest points  but  infinitely many instantons and closed trajectories.  However, 
despite their infiniteness, they can be counted by appropriate counting functions 
which can be  related to  the topology and  the geometry  of the underlying manifold cf.\ \cite{BH03}. 
The dynamical torsion is derived from them.

All torsion functions referred to above involve some additional topological data; 
the Milnor--Turaev and dynamical torsion involve an Euler structure  while the 
complex Ray--Singer torsion a coEuler structure, a sort of Poincar\'e dual of the first.  
In section~\ref{S:3} we define Euler and coEuler structures and discuss some of their 
properties. Although they can be defined for arbitrary base pointed manifolds $(M,x_0)$ 
we present the theory only in the case $\chi(M)=0$ when the base point is irrelevant.

While the complex Ray--Singer torsion and dynamical torsion are new concepts the 
Milnor--Turaev torsion is not, however our presentation is somehow different from 
the traditional one. In section~\ref{S:4} we discuss the algebraic variety 
of cochain complexes of finite dimensional vector spaces and introduce  
the Milnor torsion as  a rational function on this variety.  
The Milnor--Turaev torsion is obtained as a pull back by a characteristic map  
of this rational function.

Section~\ref{S:5} is about analytic torsions. In section~\ref{SS:5.1}, we recall  
the familiar Ray--Singer torsion slightly modified with the help of a coEuler 
structure. This is a positive real valued function defined on $\Rep^M_0(\Gamma;V)$, 
the variety of the acyclic representations. We show that this function is independent 
of the Riemannian metric, and that it is the absolute value of a rational function,
provided the coEuler structure is integral.
%%%%%%%%%%%%%%%%
%locally is the absolute value of holomorphic functions,  
%more precisely of a \emph{holomorphic closed  one m-form}.  
%but globally is not in general. However its square always is. 
%The same happens for the unmodified Ray--Singer torsion.
%%%%%%%%%%%%%%%%%%%%%
In section~\ref{SS:5.3} we introduce the complex valued Ray--Singer torsion,  
and show the relation to the first.
% its absolute value is conjecturally the square of the modified Ray--Singer torsion.
The complex Ray--Singer torsion, denoted $\mathcal S\mathcal T$, is a meromorphic 
function on a finite cover of the  space of representations and is defined analytically
using regularized determinants of elliptic operators but not self adjoint.

The  Milnor--Turaev torsion, defined in section~\ref{SS:6.1}, is associated with a smooth 
manifold, a given Euler structure and a homology orientation and is constructed  
using a smooth triangulation. Its square is  conjecturally  equal to  the 
complex Ray--Singer torsion as defined in section~\ref{SS:5.3},  
when the coEuler structure for Ray--Singer corresponds, by Poincar\'e
duality map, to the Euler structure for Milnor--Turaev. The conjecture  
is true in many relevant cases, in particular for $\dim V=1$.

Up to sign the dynamical torsion, introduced in section~\ref{SS:6.2}, is associated to 
a smooth manifold and a given Euler structure and is constructed using a smooth vector 
field in the class described in section~\ref{SS:2.3}. The sign can be  fixed with the help 
of an equivalence class of orderings of the rest points of $X$, cf.\ section~\ref{SS:6.2}. 
A priori the dynamical torsion is only a partially defined holomorphic function on 
$\Rep^M(\Gamma;V)$ and is defined using the instantons and the closed trajectories of $X$.
For a representation $\rho$ the dynamical torsion is expressed as  a series which might 
not be convergent for each $\rho$ but is certainly convergent for $\rho$ in a subset 
$U$ of $\Rep^M(\Gamma;V)$ with non-empty interior. At present this convergence was 
established only in the case of rank one representations. The existence of $U$ is 
guaranteed by the exponential growth property (EG) (cf.\ section~\ref{SS:2.3} 
for the definition) required from the vector field.

The main results, Theorems~\ref{T:1}, \ref{T:2} and \ref{T:3}, establish the 
relationship between these torsion functions, at least in the case $\dim V=1$, 
and a few other relevant cases. The  same relationship is expected to hold for 
$V$ of arbitrary dimension.

One can calculate the Milnor--Turaev torsion when $M$ has a structure of mapping torus  
of a diffeomorphism $\phi$ as the ``twisted  Lefschetz zeta function''  of the 
diffeomorphism  $\phi$, cf.\ section~\ref{SS:7.1}.
The Alexander polynomial as well as the twisted Alexander polynomials of a knot 
can also be recovered from these torsions cf.\ section~\ref{SS:7.3}.
If the vector field has no rest points but admits a closed Lyapunov cohomology class,  
cf.\ section~\ref{SS:7.2}, the dynamical torsion can be expressed in terms of closed 
trajectories only, and the dynamical zeta function of the vector field (including all 
its twisted versions) can be recovered from the dynamical torsion described here.

In section~\ref{SS:8.1} we express the phase difference of the Milnor--Turaev torsion 
of two representations in the same connected component of
$\Rep^M_0(\Gamma;V)$ in terms of the Ray--Singer torsion. This invariant is analogous to 
the Atiyah--Patodi--Singer spectral flow but has not been investigated so far.
Section~\ref{S:8} discusses some progress towards a conjecture of Novikov which came 
out from the work on dynamical torsion.

\section{Characteristic forms and vector fields}\label{S:2}

\subsection{Euler, Chern--Simons, and Mathai--Quillen form}\label{SS:2.1}

Let $M$ be smooth closed manifold of dimension $n$. 
Let $\pi:TM\to M$ denote the
tangent bundle, and $\Or_M$ the orientation bundle, which is a flat real line
bundle over $M$. 
For a Riemannian metric $g$ denote by 
$$
\ec(g)\in\Omega^n(M;\Or_M)
$$
its Euler form, and for two Riemannian metrics $g_1$ and $g_2$ by 
$$
\cs(g_1,g_2)\in
\Omega^{n-1}(M;\Or_M)/d(\Omega^{n-2}(M;\Or_M))
$$ 
their Chern--Simons class. 
%The definition of both quantities is
%implicit in the formulas \eqref{dpsie} and \eqref{pg1pg2cs} below.
The following properties follow from
\eqref{dpsie} and \eqref{pg1pg2cs} below.
\begin{eqnarray}
d \cs(g_1,g_2) &=& \ec(g_2)-\ec(g_1)
\label{E:csg:i}
\\
\cs(g_2,g_1) &=& -\cs(g_1,g_2)
\label{E:csg:ii}
\\
\cs(g_1,g_3) &=& \cs(g_1,g_2) + \cs(g_2,g_3)
\label{E:csg:iii}
\end{eqnarray}
If the dimension of $M$ is odd both $\ec(g)$ and $\cs(g_1, g_2)$ vanish.

Denote by $\xi$ the Euler vector field on $TM$ which
assigns to a point $x\in TM$ the vertical vector $-x\in T_xTM$.
A Riemannian metric $g$ determines the Levi--Civita connection in the bundle
$\pi:TM\to M$. There is a canonic $n$-form $\Vol(g)\in\Omega^n(TM;\pi^*\Or_M)$,
which assigns to an $n$-tuple of vertical vectors \emph{their volume times their
orientation} and vanishes when contracted with horizontal vectors and a global 
angular form, see for instance \cite{BT82}, is the differential form
$$
A(g)
:=\frac{\Gamma(n/2)}{(2\pi)^{n/2}|\xi|^n}
i_\xi\Vol(g)\in\Omega^{n-1}(TM\setminus M;\pi^*\Or_M).
$$
In \cite{MQ86} Mathai and Quillen have introduced a differential form 
$$
\Psi(g)\in \Omega^{n-1}(TM\setminus M;\pi^*\Or_M).
$$ 
When pulled back to the fibers of $TM\setminus M\to M$ the form $\Psi(g)$
coincides with $A(g)$. If $U\subseteq M$ is an open subset on
which the curvature of $g$ vanishes, then $\Psi(g)$ coincides with $A(g)$
on $TU\setminus U$.
In general we have the equalities
\begin{eqnarray}
d\Psi(g) &=& \pi^*\ec(g).
\label{dpsie}
\\
\Psi(g_2)-\Psi(g_1) &=& \pi^*\cs(g_1,g_2) \mod d\Omega^{n-2}(TM \setminus M;\pi^* \Or_M).
\label{pg1pg2cs}
\end{eqnarray}

\subsection{Euler and Chern--Simons chains}\label{SS:2.22}

For a vector field $X$ with non-degenerate rest points 
we have the singular $0$-chain $\ec(X)\in C_0(M;\mathbb Z)$ defined by 
$\ec(X):=\sum_{x\in\mathcal X}\IND(x)x$, with $\IND(x)$ the Hopf index.

For two  vector fields $X_1$ and $X_2$ with non-degenerate rest points 
we have the singular $1$-chain rel.\ boundaries 
$\cs(X_1,X_2)\in C_1(M;\mathbb Z)/\partial C_2(M;\mathbb Z)$ defined from the 
zero set of a homotopy from $X_1$ to $X_2$ 
cf.\ \cite{BH04}. They are related by the formulas, see \cite{BH04},
\begin{eqnarray}
\partial \cs(X_1,X_2) &=& \ec(X_2)-\ec(X_1)
\label{E:11}
\\
\cs(X_2,X_1) &=& -\cs(X_1,X_2)
\label{E:12i}
\\
\cs(X_1,X_3) &=& \cs(X_1,X_2) + \cs(X_2,X_3).
\label{E:13}
\end{eqnarray}

\subsection{Kamber--Tondeur one form}\label{SS:2.2}

Let $E$ be a real or complex vector bundle over $M$.
For a connection $\nabla$ and a Hermitian structure $\mu$ on $E$
define a real valued one form
$\omega(\nabla,\mu)\in\Omega^1(M;\R)$ by
\begin{equation}\label{E:6}
\omega(\nabla,\mu)(Y):=-\frac12\tr\bigl(\mu^{-1}\cdot (\nabla_Y\mu)\bigr),
\quad Y\in TM.
\end{equation}
Here we consider $\mu$ as an element in $\Omega^0(M;\hom(E,\bar E^*))$, where
$\bar E^*$ denotes the dual of the complex conjugate bundle. With respect to the 
induced connection on $\hom(E,\bar E^*)$ we have $\nabla_Y\mu\in\Omega^1(M;\hom(E,\bar E^*))$
and therefore $\mu^{-1}\cdot\nabla_Y\mu\in\Omega^1(M;\eend(E,E))$. Actually the latter
one form has values in the endomorphisms of $E$ which are symmetric with respect to $\mu$,
and thus the (complex) trace, see \eqref{E:6}, will indeed be real.
Since any two Hermitian structures $\mu_1$ and $\mu_2$ are homotopic,
the difference $\omega(\nabla,\mu_2)-\omega(\nabla,\mu_1)$ will be exact.
If $\nabla$ is flat then $\omega(\nabla,\mu)$ is 
closed and its cohomology class independent of $\mu$.
%Moreover, 
%$\omega(\nabla^{\det E},\mu^{\det E})=\omega(\nabla,\mu)$
%where $\nabla^{\det E}$ and $\mu^{\det E}$ denote the induced connection
%and Hermitian structure on the determinant line 
%$\det E:=\Lambda^{\rank E}E$.

%A similar form (closed for a flat connections) can be associated to a symmetric
%non-degenerate bilinear form $b$
%on $E$ instead of a Hermitian structure by a similar formula
Replacing the Hermitian structure by a non-degenerate symmetric bilinear form $b$,
we define a complex valued one form $\omega(\nabla,b)\in\Omega^1(M;\C)$ by
a similar formula
\begin{equation}\label {E:7}
\omega(\nabla,b)(Y):=-\frac12\tr\bigl(b^{-1}\cdot(\nabla_Yb)\bigr),\quad Y\in TM.
\end{equation}
Here we regard $b$ as an element in $\Omega^0(M;\hom(E,E^*))$.
If two non-degenerate symmetric bilinear forms $b_1$ and $b_2$ are homotopic, then 
$\omega(\nabla,b_2)-\omega(\nabla,b_1)$ is exact. If $\nabla$ is flat,
then $\omega(\nabla,b)$ is closed.
Note that $\omega(\nabla,b)\in\Omega^1(M;\mathbb C)$ depends holomorphically
on $\nabla$.

\subsection{Vector fields, instantons and closed trajectories}\label{SS:2.3}

Consider a vector field $X$ which satisfies the following properties:

\begin{itemize}
\item[(H)]
All rest points are of hyperbolic type.
\item[(EG)]
The vector field has exponential growth at all rest points.
\item[(L)]
The vector field is of Lyapunov type.
\item[(MS)]
The vector field satisfies Morse--Smale condition.
\item[(NCT)]
The vector field has all closed trajectories non-degenerate.
\end{itemize}
Precisely this means that:
\begin{itemize}
\item[(H)]
In the neighborhood of each rest point the differential of $X$ has all 
eigenvalues with non-trivial real part; the number of eigenvalues with 
negative real part is called the index and denoted by $\ind(x)$; 
as a consequence the stable
and unstable stable sets are images of one-to-one immersions 
$i_x^\pm:W^\pm_x\to M$ with $W^\pm_x$ diffeomorphic to $\mathbb R^{n-\ind(x)}$ 
resp.\ $\mathbb R^{\ind(x)}$.
\item[(EG)]
With respect to one and then any Riemannian metric $g$ on
$M$, the volume of the disk of radius $r$ in $W^-_x$ (w.r.\ to the induced Riemannian
metric) is $\leq e^{Cr}$, for some constant $C>0$.
\item[(L)]
There exists a real valued closed one form $\omega$ so that $\omega(X)_x<0$ for 
$x$ not a rest point.\footnote{This $\omega$ has nothing in common with 
$\omega(\nabla,b)$ notation used in the previous section.}
\item[(MS)]
For any two rest points $x$ and $y$ the mappings $i^-_x$ and $i^+_y$ are transversal 
and therefore the space of non-parameterized trajectories form $x$ to $y$, 
$\mathcal T(x,y)$, is a smooth manifold of dimension $\ind(x)-\ind(y)-1$. 
Instantons are exactly the elements of $\mathcal T(x,y)$ when this is a smooth 
manifold of dimension zero, i.e.\ $\ind(x)-\ind(y)-1=0$.
\item[(NCT)]
Any closed trajectory is non-degenerate, i.e.\ the differential of the return map
in normal direction at one and then any point of a closed trajectory does not have
non-zero fixed points.
\end{itemize}

Recall that a trajectory $\theta$ is an equivalence class of parameterized trajectories and 
two parameterized trajectories $\theta_1$ and $\theta_2$ are equivalent iff 
  $\theta_1(t+c)= \theta_2(t)$
for some real number $c.$  
Recall that a closed trajectory $\hat\theta$ is a pair consisting of a trajectory $\theta$ and a positive real number $T$ so that $\theta(t+T)= \theta(t).$ 

%an equivalence class of parameterized closed trajectory 
%consisting of ais a pair $(\theta, T)$ with $\theta$ a parameterized 
%trajectory  and $T$ a positive real number so that $\theta(t)= \theta(t+T).$  A closed trajectory $\hat %\theta$ is an equivalence class and parameterized closed trajectories;  $(\theta_1, T_1)\equiv %(\theta_2, T_2)$ iff $T_1= T_2$ and $\theta_1$ and $\theta_2$ are equivalent parameterized %trajectories. 

Property (L), (H), (MS)  imply that for any real number $R$ the set of instantons 
$\theta$ from $x$ to $y$ 
with $-\omega([\theta])\leq R$ 
is finite  and properties (L), (H), (MS), (NCT) imply that for any real number $R$ the set of 
the closed trajectory $\hat\theta$ with $-\omega([\hat \theta])\leq R$ is finite. 
Here 
$[\theta]$ resp $[\hat\theta]$ denote the homotopy class of instantons resp. closed trajectories
\footnote{ For a closed trajectory the map whose homotopy class is considered is $\hat \theta: \mathbb R/ T\mathbb Z\to M.$}.

Denote by $\mathcal P _{x,y}$ the set of homotopy classes of paths from 
$x$ to $y$ and by $\mathcal X_q$ the set of rest points of index $q$. 
Suppose a collection $\mathcal O=\{\mathcal O_x\mid x\in\mathcal X\}$
of orientations of the unstable manifolds is given and (MS) is satisfied.
Then any instanton $\theta$ has a sign $\epsilon(\theta)=\pm1$ and therefore,
if (L) is also satisfied, for any two rest points $x\in{\mathcal X}_{q+1}$
and $y\in \mathcal X_{q}$ we have the counting function of instantons
$\mathbb I^{X,\mathcal O}_{x,y}: \mathcal P_{x,y}\to\mathbb Z$
defined by 
\begin{equation}\label {E:8}
\mathbb I^{X,\mathcal O}_{x,y}(\alpha):=\sum_{\theta\in\alpha}\epsilon(\theta).
\end{equation}

Under the hypothesis (NCT) any closed trajectory $\hat\theta$ has a sign 
$\epsilon(\hat\theta)=\pm1$ and a period $p(\hat\theta)\in\{1,2,\dotsc\}$, cf.\ \cite{H02}.
If (H), (L), (MS), (NCT) are satisfied, as the set of closed trajectories in a fixed 
homotopy class $\gamma\in[S^1,M]$ is compact, we have the counting function of closed 
trajectories $\mathbb Z_X:[S^1,M]\to\mathbb Q$ defined by  
\begin{equation}\label{E:9}
\mathbb Z_X(\gamma):=\sum_{\hat\theta\in\gamma}\epsilon(\hat\theta)/p(\hat\theta).
\end{equation}
Here are a few properties about vector fields which satisfy (H) and (L).

\begin{proposition}\label{P:1}
1. Given a vector field $X$ which satisfies (H) and (L) arbitrary close in the
$C^r$-topology for any $r\geq0$ there exists a vector field $Y$ which agrees  
with $X$ on a neighborhood of the rest points and satisfies (H), (L), (MS) and (NCT).

2. Given a vector field $X$ which satisfies (H) and (L) arbitrary close in the
$C^0$-topology there exists a vector field $Y$ which agrees with $X$ on a 
neighborhood of the rest points and satisfies (H), (EG), (L), (MS) and (NCT).

3. If $X$ satisfies (H), (L) and (MS) and a collection $\mathcal O$ of orientations 
is given then  for any  $x\in\mathcal X_q$, $z\in\mathcal X_{q-2}$ and 
$ \gamma\in\mathcal P_{x,z}$ one has\footnote{It is understood that only finitely many 
terms from the left side of the equality are not zero. Here $*$ denotes juxtaposition.}
\begin{equation}\label{E:10}
%\sum_{ \begin{aligned} y\in{\mathcal X}_{q-1},   \alpha\in &\mathcal P_{x,y},  
%\beta\in \mathcal P_{y,z} \\  \alpha* \beta&= \gamma\end{aligned}}
\sum_{\begin{smallmatrix} y\in{\mathcal X}_{q-1},\alpha\in\mathcal P_{x,y},  
\beta\in \mathcal P_{y,z} \\  \alpha* \beta=\gamma\end{smallmatrix}}
\mathbb I^{X,\mathcal O}_{x,y}(\alpha)\cdot
\mathbb I^{X,\mathcal O}_{y,z}(\beta)=0.
\end{equation}
\end{proposition}

This proposition is a recollection of some of the main results in \cite{BH03}, see
Proposition~3, Theorem~1 and Theorem~5 in there.

\section{Euler and coEuler structures}\label{S:3}

Although not always necessary in this section as in fact always in this 
paper $M$ is supposed to be closed connected smooth manifold.

\subsection{Euler structures}\label{SS:3.1}

Euler structures have been introduced by Turaev \cite{Tu90} for
manifolds $M$ with $\chi(M)=0$. If one
removes the hypothesis $\chi(M)=0$ the concept of Euler structure 
can still be considered for any connected base pointed manifold
$(M,x_0)$ cf.\ \cite{B99} and \cite{BH04}. Here we will consider only
the case $\chi(M)=0$. The set of Euler structures,
denoted by $\Eul(M;\Z)$, is equipped with a free and transitive action
$$
m:H_1(M;\Z)\times\Eul(M;\Z)\to\Eul(M;\Z) 
$$ 
which makes $\Eul(M;\Z)$ an affine version of $H_1(M;\Z)$.
If $\e_1,\e_2$ are two Euler structure we write $\e_2-\e_1$ for the 
unique element in $H_1(M;\Z)$ with $m(\e_2-\e_1,\e_1)=\e_2$.

To define the set $\Eul(M;\Z)$ we consider pairs $(X,c)$ with $X$ a vector field with
non-degenerate zeros and $c\in C_1(M;\mathbb Z)$ so that $\partial
c=\ec(X)$. We make $(X_1,c_1)$ and $(X_2,c_2)$ equivalent iff $c_2-c_1=\cs(X_1,X_2)$
and write $[X,c]$ for the equivalence class represented by $(X,c)$.
The action $m$ is defined by $m([c'],[X, c]):=[X,c'+c]$.

\begin{observation}\label{O:1}
Suppose $X$ is a vector field with non-degenerate zeros, and assume its zero set
$\mathcal X$ is non-empty. Moreover, let $\e\in\Eul(M;\mathbb Z)$ be an Euler structure
and $x_0\in M$. Then there exists a collection of
paths $\{\sigma_x\mid x\in\mathcal X\}$ with $\sigma_x(0)=x_0$,
$\sigma_x(1)=x$ and such that $\e=[X,c]$ where $c=\sum_{x\in\mathcal X}\IND(x)\sigma_x$.
\end{observation}

A remarkable source of Euler structures is the set of homotopy classes
of nowhere vanishing vector fields. Any nowhere vanishing vector field $X$
provides an Euler structure $[X,0]$ which only depends on the homotopy class
of $X$. Still assuming $\chi(M)=0$, every 
Euler structure can be obtained in this way provided $\dim(M)>2$.
Be aware, however, that different homotopy classes may give rise to the
same Euler structure.

To construct such a homotopy class one can proceed as follows. Represent
the Euler structure $\e$ by a vector field $X$ and a collection of paths
$\{\sigma_x\mid x\in\mathcal X\}$ as in Observation~\ref{O:1}.
Since $\dim(M)>2$ we may assume that the interiors
of the paths are mutually disjoint. Then the set $\bigcup_{x\in\mathcal
X}\sigma_x$ is contractible. A smooth regular neighborhood
of it is the image by a smooth embedding $\varphi:(D^n,0)\to(M,x_0)$. 
Since $\chi(M)=0$, the restriction of the vector field $X$ to
$M\setminus\inter(D^n)$ can be extended to a non-vanishing vector field 
$\tilde X$ on $M$. It is readily checked that $[\tilde X,0]=\e$.
For details see \cite{BH04}.

If $M$  dimension larger than 
2 an alternative description of $\Eul(M;\Z)$
%%%%%%%%%%%%%%
% and of the action $m$  
%%%%%%%%%%%
with respect to a base point $x_0$ is 
$\Eul(M;\Z)=\pi_0(\vf(M,x_0))$, where
$\vf(M,x_0)$ denotes the space of vector fields of class $C^r$, $r\geq 0$, which vanish 
at $x_0$ and are non-zero elsewhere. We equip this space with the $C^r$-topology and 
note that the result $\pi_0(\vf(M,x_0))$ is the same for all $r$, and
since $\chi(M)=0$, canonically identified for different base points.

Let $\tau$ be a smooth triangulation of $M$ and consider the function 
$f_\tau:M\to\R$ linear on any simplex of the first barycentric subdivision 
and taking the value $\dim(s)$ on the barycenter $x_s$ of
the simplex $s\in\tau$. A smooth vector field $X$ on $M$ 
with the barycenters as the only rest points all of them hyperbolic 
%as hyperbolic rest points
and $f_\tau$ strictly decreasing on
non-constant trajectories is called an Euler vector field of $\tau$.
By an argument of convexity two Euler vector fields are homotopic by a homotopy of
Euler vector fields.\footnote{Any Euler vector field $X$ satisfies (H), (EG), (L) 
and has no closed trajectory, hence also satisfies (NCT).
The counting functions of instantons are exactly the same as the incidence numbers  
of the triangulation hence take the values $1$, $-1$ or $0$.}
Therefore, a triangulation $\tau$, a base point $x_0$ and a collection of
paths $\{\sigma_s\mid s\in\tau\}$ with $\sigma_s(0)=x_0$ and
$\sigma_s(1)=x_s$ define an Euler structure $[X_\tau,c]$, where
$c:=\sum_{s\in\tau}(-1)^{n+\dim(s)}\sigma_s$, $X_\tau$ is any Euler vector
field for $\tau$, and this Euler structure does not depend on the
choice of $X_\tau$. Clearly, for fixed $\tau$ and $x_0$, every Euler
structure can be realized in this way by an appropriate choice of
$\{\sigma_s\mid s\in\tau\}$, cf.\ Observation~\ref{O:1}.

\subsection{Co-Euler structures}\label{SS:3.2}

Again, suppose $\chi(M)=0$.\footnote{The hypothesis is not necessary and the theory 
of coEuler structure can be pursued for an arbitrary base pointed smooth manifold
$(M,x_0)$, cf.\ \cite{BH04}.} Consider pairs
$(g,\alpha)$ where $g$ is a Riemannian metric on $M$ and
$\alpha\in\Omega^{n-1}(M;\mathcal O_M)$ with $d\alpha=\ec(g)$
where $\ec(g)\in\Omega^n(M;\mathcal O_M)$ denotes the Euler form of 
$g$, see section~\ref{SS:2.1}. We call two
pairs $(g_1,\alpha_1)$ and $(g_2,\alpha_2)$ equivalent if
$$
\cs(g_1,g_2)=\alpha_2-\alpha_1
\in\Omega^{n-1}(M;\mathcal O_M)/
d\Omega^{n-2}(M;\mathcal O_M).
$$
We will write $\Eul^*(M;\R)$ for the set of equivalence classes and
$[g,\alpha]$ for the equivalence class represented by the pair $(g,\alpha)$.
Elements of $\Eul^*(M;\R)$ are called \emph{coEuler structures.}

There is a natural action
$$
m^*:H^{n-1}(M;\Or_M)\times\Eul^*(M;\R)\to\Eul^*(M;\R)
$$
given by
$$
m^*([\beta],[g,\alpha]):=[g,\alpha-\beta]
$$
for $[\beta]\in H^{n-1}(M;\mathcal O_M)$. 
This action is obviously free and transitive. In this sense
$\Eul^*(M;\R)$ is an affine version of $H^{n-1}(M;\mathcal O_M)$.
If $\e^*_1$ and $\e^*_2$ are two coEuler structures we write $\e^*_2-\e^*_1$
for the unique element in $H^{n-1}(M;\Or_M)$ with 
$m^*(\e^*_2-\e^*_1,\e^*_1)=\e^*_2$.

\begin{observation}\label {O:2}
Given a Riemannian metric $g$ on $M$ any coEuler structure can be
represented as a pair $(g,\alpha)$ 
for some $\alpha\in\Omega^{n-1}(M;\Or_M)$ with $d\alpha=\ec(g)$.
\end{observation}

There is a natural map $\PD:\Eul(M;\Z)\to\Eul^*(M;\R)$
which combined with the Poincar\'e duality map
$D:H_1(M;\Z)\to H_1(M;\R)\to H^{n-1}(M;\Or_M)$, the composition of the coefficient
homomorphism for $\Z\to\R$ with the Poincar\'e duality isomorphism,\footnote{We 
will use the same notation D for the Poincar\'e duality isomorphism
$D:H_1(M;\R)\to H^{n-1}(M; \mathcal O_M)$.} makes the diagram below commutative:
$$
\xymatrix{
H_1(M;\Z)\times\Eul(M;\Z)
\ar[d]_{D\times\PD}\ar[r]^-{m}   
& 
\Eul(M;\Z)\ar[d]^{\PD}
\\
H^{n-1}(M;\Or_M)\times\Eul^*(M;\R) 
\ar[r]^-{m^*} 
&  
\Eul^*(M;\R)
}
$$

There are many ways to define the map $\PD$, cf.\ \cite{BH04}. For example, 
assuming $\chi(M)=0$ and $\dim M>2$ one can proceed as follows.
Represent the Euler structure by a nowhere vanishing vector field
$\e=[X,0]$. Choose a Riemannian metric $g$, regard $X$ as 
mapping $X:M\to TM\setminus M$, set $\alpha:=X^*\Psi(g)$, put 
$\PD(\e):=[g,\alpha]$ and check that this does indeed only depend on $\e$.

A coEuler structure $\e^*\in\Eul^*(M;\R)$ is called \emph{integral}
if it belongs to the image of $\PD$. Integral coEuler structures 
constitute a lattice in the affine space $\Eul^*(M;\R)$.

\begin{observation} If $\dim M$ is odd, then there is a canonical coEuler structure
$\e^*_0\in\Eul^*(M;\R)$; it is represented by the pair $[g,0]$, with any $g$ 
Riemannian metric. In general this coEuler structure is not integral. 
\end{observation}

\section{Complex representations and cochain complexes}\label{S:4}

\subsection{Complex representations}\label{SS:4.1}

Let $\Gamma$ be a finitely presented group with generators $g_1,\dotsc,g_r$ 
and relations 
$$
R_i(g_1,g_2,\dotsc,g_r)=e,\quad i=1,\dotsc,p,
$$ 
and $V$ be a
complex vector space of dimension $N$. Let $\Rep(\Gamma;V)$ be the set of 
linear representations of $\Gamma$ on $V$, i.e.\ group homomorphisms 
$\rho:\Gamma\to\GL_\C(V)$. By identifying $V$ to $\C^N$ this set is, in a
natural way, an algebraic set inside the space $\C^{rN^2+1}$ given
by $pN^2+1$ equations. Precisely if $A_1,\dotsc,A_r,z$ represent the 
coordinates in $\C^{rN^2+1}$ with $A:=(a^{ij})$, $a^{ij}\in\C$, so 
$A\in\C^{N^2}$ and $z\in\C$, then the equations defining $\Rep(\Gamma;V)$
are 
\begin{eqnarray*}
z\cdot\det(A_1)\cdot\det(A_2)\cdots\det(A_r)&=&1
\\
R_i(A_1,\dotsc,A_r)&=&\id,\qquad i=1,\dotsc,p
\end{eqnarray*}
with each of the equalities $R_i$ representing $N^2$ polynomial equations.

Suppose $\Gamma=\pi_1(M,x_0)$, $M$ a closed manifold.
Denote by $\Rep^M_0(\Gamma;V)$ the set of representations $\rho$ with 
$H^*(M;\rho)=0$ and notice that they form a Zariski
open set in $\Rep(\Gamma;V)$. Denote the closure of this set by  
$\Rep^M(\Gamma;V)$. This is an algebraic set 
which depends only on the homotopy type of $M$,
and is a union of irreducible components of $\Rep(\Gamma;V)$.

Recall that every representation $\rho\in\Rep(\Gamma;V)$ induces a
canonical vector bundle $F_\rho$ equipped with a canonical flat
connection $\nabla_\rho$. They are obtained from the
trivial bundle $\tilde M\times V\to\tilde M$ and the trivial connection by
passing to the $\Gamma$ quotient spaces. Here $\tilde M$ is the canonical
universal covering provided by the base point $x_0$. The $\Gamma$-action is
the diagonal action of deck transformations on $\tilde M$ and of the action 
$\rho$ on $V$. The fiber of $F_\rho$ over $x_0$ identifies canonically with
$V$. The holonomy representation determines a right $\Gamma$-action on the
fiber of $F_\rho$ over $x_0$, i.e.\ an anti homomorphism $\Gamma\to\GL(V)$.
When composed with the inversion in $\GL(V)$ we get back the representation
$\rho$. The pair $(F_\rho,\nabla_\rho)$ will be denoted by $\mathbb F_\rho$.

If $\rho_0$ is a representation in the connected component 
$\Rep_\alpha(\Gamma;V)$ one can identify $\Rep_\alpha(\Gamma;V)$ to the
connected component of $\nabla_{\rho_0}$ in the complex analytic space of
flat connections of the bundle $F_{\rho_0}$ modulo the group of bundle
isomorphisms of $F_{\rho_0}$ which fix the fiber above $x_0$.

\begin{remark}
An element $a\in H_1(M;\Z)$ defines a holomorphic function  
$$
{\det}_a:\Rep^M(\Gamma;V)\to\mathbb C_*.
$$
The complex number $\det_a(\rho)$ is
the evaluation on $a\in H_1(M;\Z)$ of $\det(\rho):\Gamma\to\C_*$
which factors through $H_1(M;\Z)$. Note that for $a,b\in H_1(M;\Z)$ we have 
$\det_{a+b}=\det_a\det_b$. 
If $a$ is a torsion element, then $\det_a$ is constant equal
to a root of unity of order, the order of $a$. 
\end{remark}

\subsection{The space of cochain complexes}\label{SS:4.3}

Let $k= (k_0,k_1,\dotsc,k_n)$ be a string of non-negative integers. The string is 
called admissible, and will write $k\geq0$ in this case, 
if the following requirements are satisfied
\begin{eqnarray}
k_0-k_1+k_2\mp\cdots+(-1)^nk_n&=&0
\label{E:14}
\\
k_i-k_{i-1}+ k_{i-2}\mp\cdots+(-1)^ik_0&\geq&0
\quad\text{for any $i\leq n-1$.}
\label{E:15}
\end{eqnarray}

Denote by $\mathbb D(k)=\mathbb D(k_0,\dotsc,k_n)$
the collection of cochain complexes of the form
$$
C=(C^*,d^*):
0\to C^0 \xrightarrow{d^0} 
C^1\xrightarrow{d^1}
\cdots 
\xrightarrow{d^{n-2}}
C^{n-1}
\xrightarrow{d^{n-1}} C^n\to 0
$$
with $C^i:=\C^{k_i}$, and by $\mathbb D_\ac(k)\subseteq\mathbb D(k)$ 
the subset of acyclic complexes. Note that $\mathbb D_\ac(k)$ is non-empty
iff $k\geq0$. The cochain complex $C$ is determined by 
the collection $\{d^i\}$ of linear maps $d^i:\C^{k_i}\to\C^{k_{i+1}}$.
If regarded as the subset of those
$\{d^i\}\in\bigoplus_{i=0}^{n-1}L(\C^{k_i},\C^{k_{i+1}})$, 
with $L(V,W)$ the space of linear maps from $V$ to $W$,
which satisfy the quadratic equations $d^{i+1}\cdot d^i=0$,
the set $\mathbb D(k)$ is an affine algebraic 
set given by degree two homogeneous polynomials and 
$\mathbb D_\ac(k)$ is a Zariski open set.
The map $\pi_0:\mathbb D_\ac(k)\to\Emb(C^0,C^1)$
which associates to $C\in\mathbb D_\ac(k)$
the linear map $d^0$, is a bundle whose fiber is isomorphic to 
$\mathbb D_\ac(k_1-k_0,k_2,\dotsc,k_n)$.

This can be easily generalized as follows.
%%%%%%%%%%%
%Consider two $n$-strings $k$ and $b=(b_0,\dotsc,b_n)$ and suppose that 
%$k-b=(k_0-b_0,\dotsc,k_n-b_n)$ is admissible. We write this as $k\geq b$.  
%%%%%%%%%%%
Consider a string $b=(b_0,\dotsc,b_n)$. We will write $k\geq b$ if
$k-b=(k_0-b_0,\dotsc,k_n-b_n)$ is admissible, i.e.\ $k-b\geq0$.
Denote by $\mathbb D_b(k)=\mathbb D_{(b_0,\dotsc,b_n)}(k_0,\dotsc,k_n)$ 
the subset of cochain complexes $C\in\mathbb D(k)$ with $\dim(H^i(C))=b_i$.
Note that $\mathbb D_b(k)$ is non-empty iff $k\geq b$.
The obvious map $\pi_0:\mathbb D_b(k)\to L(C^0,C^1;b_0)$, 
$L(C^0,C^1;b_0)$ the space of linear maps in $L(C^0,C^1)$ whose kernel 
has dimension $b_0$, 
%%%%%%%%%%%%%%%%%%%%%%%%%%%
%%$ 
%\pi_0:\mathbb D_\ac (k_0,\dotsc, k_n)\to   \Emb (\ker d^0, C^0)\times \Emb(C^0/ \ker d^0 ,C^1)
%$
%%%%%%%%%%%%%%%%%%%%%%%%%%%%%%
is a bundle whose fiber is isomorphic to 
$\mathbb D_{b_1,\dotsc,b_n}(k_1-k_0+b_0,k_2,\dotsc,k_n)$. Note that
$L(C^0, C^1; b_0)$ is the total space of a bundle
$\Emb(\underline{\C}^{k_0}/L,\underline{\C}^{k_1})\to\Gr_{b_0}(k_0)$ with
$L\to\Gr_{b_0}(k_0)$ the tautological bundle over $\Gr_{b_0}(k_0)$ and 
$\underline{\C}^{k_0}$ resp.\ $\underline{\C}^{k_1}$ the trivial bundles over 
$\Gr_{b_0}(k_0)$ with fibers of dimension $k_0$ resp.\ $k_1$.
As a consequence we have
\begin{proposition}\label{P:2}
1. $\mathbb D_\ac(k)$ and $\mathbb D_b(k)$ are connected smooth quasi affine  
algebraic sets whose dimension is 
$$
\dim\mathbb D_b(k)=\sum_j(k^j-b^j)\cdot\Bigl(
k^j-\sum_{i\leq j}(-1)^{i+j}(k^i-b^i)
\Bigr).
$$

2. The closures $\hat{\mathbb D}_\ac(k)$ and $\hat{\mathbb D}_b(k)$ are irreducible 
algebraic sets, hence affine algebraic varieties, and
$\hat{\mathbb D}_b(k)=\bigsqcup_{k\geq b'\geq b}{\mathbb D}_{b'}(k)$.
\end{proposition}

For any cochain complex in $C\in\mathbb D_\ac(k)$ denote by 
$B^i:=\img(d^{i-1})\subseteq C^i= \C^{k_i}$
and consider the short exact sequence 
$0\to B^i \xrightarrow{\textrm{inc}} C^i \xrightarrow {d^i} B^{i+1}\to 0$. 
Choose a base $\bb_i$ for each $B_i$, and choose lifts $\overline{\bb}_{i+1}$ of
$\bb_{i+1}$ in $C^i$ using $d^i$, i.e.\ $d^i(\overline{\bb}_{i+1})=\bb_{i+1}$. 
Clearly $\{\bb_i,\overline\bb_{i+1}\}$ is a base of $C^i$.
Consider the base $\{\bb_i,\overline{\bb}_{i+1}\}$ as a collection of vectors in 
$C^i=\C^{k_i}$ and write them as columns of a matrix 
$[\bb_i,\overline{\bb}_{i+1}]$. Define the torsion of the acyclic complex $C$, by 
$$
\tau(C):=(-1)^{N+1}\prod^n_{i=0}\det[\bb_i,\overline{\bb}_{i+1}]^{(-1)^i}
$$
where $(-1)^N$ is Turaev's sign, see \cite{FT99}.
The result is independent of the choice of 
the bases $\bb_i$ and of the lifts $\overline{\bb}_i$ cf.\ \cite{M66}
\cite{FT99}, and leads to the function 
$$
\tau:\mathbb D_\ac(k)\to\C_*.
$$ 
Turaev provided a simple formula for this function, cf.\ \cite{Tu01}, which  
permits to recognize $\tau$ as the restriction of a rational function on  
$\hat{\mathbb D}_\ac(k)$.

For $C\in\hat{\mathbb D}_\ac(k)$ denote by $(d^i)^t:\C^{k_{i+1}}\to\C^{k_i}$ 
the transpose of $d^i:\C^{k_i}\to\C^{k_{i+1}}$, and define 
$P_i=d^{i-1}\cdot(d^{i-1})^t+(d^i)^t\cdot d^i$. Define $\Sigma(k)$ as the 
subset of cochain complexes in $\hat{\mathbb D}_\ac(k)$ where $\ker P\neq0$, and consider
$S\tau:\hat{\mathbb D}_\ac(k)\setminus\Sigma(k)\to\mathbb C_*$ defined by 
$$
S\tau(C)
:=\Bigl(\prod_{i\ \even}(\det P_i)^{i}
\big/\prod_{i\ \odd}(\det P_i)^{i}\Bigr)^{-1}.
$$
One can verify

\begin{proposition} 
Suppose $k=(k_0,\dotsc,k_n)$ is admissible.

1. $\Sigma(k)$ is a proper subvariety containing the singular set of 
$\hat{\mathbb D}_\ac(k)$.

2. $S\tau=\tau^2$ and implicitly $S\tau$ has an analytic continuation 
to ${\mathbb D}_\ac(k)$.
\end{proposition}

In particular $\tau$ defines a square root of $S\tau$.  
We will not use explicitly $S\tau $ in this writing
however it justifies the definition of complex Ray--Singer torsion.

\section{Analytic torsion}\label{S:5}

Let $M$  be a closed manifold, $g$ Riemannian metric and $(g,\alpha)$ 
a representative of a coEuler structure $\e^*\in\Eul^*(M;\R)$.
Suppose $E\to M$ is a complex vector bundle and denote by $\mathcal C(E)$ 
the space of connections and by $\mathcal F(E)$ the subset of flat connections.
$\mathcal C(E)$ is a complex affine (Fr\'echet) space while $\mathcal F(E)$  
a closed complex analytic subset (Stein space) of $\mathcal C(E)$.
Let $b$ be a non-degenerate symmetric bilinear form and $\mu$ a Hermitian 
(fiber metric) structure on $E$.  While Hermitian structures always exist,  
non-degenerate symmetric bilinear forms exist 
iff the bundle is the complexification of some real vector bundle, and in
this case $E\simeq E^*$.

The connection $\nabla\in\mathcal C(E)$ can be interpreted as a first order 
differential operator
$d^\nabla:\Omega^*(M;E)\to\Omega^{*+1}(M;E)$ and $g$ and $b$ resp.\ $g$ and $\mu$
can be used to define the formal $b$-adjoint resp.\ $\mu$-adjoint 
$\delta^\nabla_{q;g,b}$ resp.\ $\delta^\nabla_{q;g,\mu}:\Omega^{q+1}(M;E)\to\Omega^q(M;E)$
and therefore the Laplacians
$$
\Delta^\nabla_{q;g,b} 
\ \text{resp.}\
\Delta^\nabla_{q;g,\mu}:\Omega^q(M;E)\to\Omega^q(M;E).
$$
They are elliptic second order differential operators with principal symbol 
$\sigma_\xi=|\xi|^2$. Therefore they have a
unique well defined zeta regularized determinant (modified determinant) 
$\det(\Delta^\nabla_{q;g,b})\in\C$ (${\det}'(\Delta^\nabla_{q;g,b})\in\C_*$)
resp.\ $\det(\Delta^\nabla_{q;g,\mu})\in\R_{\geq 0}$  
(${\det}'(\Delta^\nabla_{q;g,\mu})\in\R_{>0}$) calculated with respect to
a non-zero Agmon angle avoiding the spectrum cf.\ \cite{BH05}.
% Agmon angle $\pi$.
%%%%%%%%%%%%
%This determinant is zero if $\ker\Delta^\nabla_q\neq 0$.
%%%%%%%%%%%%%%%
Recall that the zeta regularized determinant (modified determinant) is the 
zeta regularized product of all (non-zero) eigenvalues.
%%%%%%%%%%%%%%
%(sometimes called  the modified zeta regularized determinant.)
%%%%%%%%%%%%%%%

Denote by 
\begin{align*}
\Sigma(E,g,b)
&:=\bigl\{\nabla\in\mathcal C(E)\bigm| 
\ker(\Delta^\nabla_{*;g,b})\neq0\bigr\}
\\
\Sigma(E,g,\mu)
&:=\bigl\{\nabla\in\mathcal C(E)\bigm| 
\ker(\Delta^\nabla_{*;g,\mu})\neq0\bigr\}
\end{align*}
and by 
$$
\Sigma(E):=\bigl\{\nabla\in\mathcal F(E) 
\bigm| H^*(\Omega^*(M;E),d^\nabla)\neq 0\bigr\}.
$$
Note that $\Sigma(E,g,\mu)\cap\mathcal F(E)=\Sigma(E)$ for any $\mu$,
and $\Sigma(E,g,b)\cap\mathcal F(E)\supseteq\Sigma(E)$. Both,
$\Sigma(E)$ and $\Sigma(E,g,b)\cap\mathcal F(E)$, are closed complex 
analytic subsets of $\mathcal F(E)$, and
$\det(\Delta^\nabla_{q;g,\cdots})
={\det}'(\Delta^\nabla_{q;g,\cdots})$ on $\mathcal F(E)\setminus \Sigma(E,g,\cdots)$.

We consider the real analytic functions:
$T^\even_{g,\mu}:\mathcal C(E)\to\mathbb R_{\geq 0}$, 
$T^\odd_{g,\mu}:\mathcal C(E)\to\mathbb R_{\geq 0}$, 
$R_{\alpha, \mu}:\mathcal C(E)\to\mathbb R_{>0}$
and the holomorphic functions 
$T^\even_{g,b}:\mathcal C(E)\to\mathbb C$,
$T^\odd_{g,b}:\mathcal C(E)\to\mathbb C$,
$R_{\alpha,b}:\mathcal C(E)\to \mathbb C_*$
defined by: 
\begin{equation}\label{E;16'}
\begin{aligned}
T^\even_{g,\cdots}(\nabla):=& \prod_{q\ \even}(\det\Delta^\nabla_{q;g,\cdots})^q,
%T'^{\,\even}_{g,\cdots}(\nabla):=\prod_{q\ \even}({\det}'\Delta^\nabla_{q;g,\cdots})^q,
\\
T^\odd_{g,\cdots}(\nabla) :=& \prod_{q\ \odd }(\det\Delta^\nabla_{q;g,\cdots})^q ,
%T'^{\,\odd}_{g,\cdots}(\nabla) := &\prod_{q\ \odd }({\det}'\Delta^\nabla_{q;g,\cdots})^q
\\
R_{\alpha,\cdots}(\nabla):= &e^{\int_M\omega(\cdots,\nabla)\wedge\alpha}.
\end{aligned}
\end{equation}
We also write $T'^{\,\even}_{g,\cdots}$ resp.\ $T'^{\,\odd}_{g,\cdots}$ for the same 
formulas with ${\det}'$ instead of $\det$. These functions are discontinuous on 
$\Sigma(E,g,\cdots)$ and coincide with $T^\even_{g,\cdots}$ resp.\ 
$T^\odd_{g,\cdots}$ on $\mathcal F(E)\setminus\Sigma(E,g,\cdots)$.
Here $\cdots$ stands for either $b$ or $\mu$.
For the  definition of real or complex analytic space/set, holomorphic/meromorphic function/map 
in infinite dimension the reader can consult \cite {D} and \cite{KM}, although the definitions used here are rather straightforward.

Let $E_r\to M$ be a smooth real vector bundle equipped with a non-degenerate symmetric
positive definite bilinear form $b_r$.
Let $\mathcal C(E_r)$ resp.\ $\mathcal F(E_r)$ the space of connections resp.\ 
flat connections in $E_r$. Denote by $E\to M$ the complexification of $E_r$, 
$E=E_r\otimes\mathbb C$, and by $b$ resp.\ $\mu$ the complexification of $b_r$ 
resp.\ the Hermitian structure  extension of $b_r$. We continue to denote by   
$\mathcal C(E_r)$ resp.\ $\mathcal F(E_r)$ the subspace of $\mathcal C(E)$ 
resp.\ $\mathcal F(E)$ consisting of connections which are complexification of connections  
resp.\ flat connections in $E_r$, and by $\nabla$ the complexification of 
the connection $\nabla\in\mathcal C(E_r)$. If $\nabla\in\mathcal C(E_r)$, then
$$
\Spect\Delta^\nabla_{q;g,b}=\Spect\Delta^\nabla_{q;g,\mu}\subseteq\mathbb R_{\geq 0}
$$
and therefore  
\begin{equation}\label{E:16}
\begin{aligned}
T^{\even/\odd}_{g,b}(\nabla)&=\bigl|T^{\even/\odd}_{g,b}(\nabla)\bigr|
=T^{\even/\odd}_{g,\mu}(\nabla),
\\
T'^{\,\even/\odd}_{g,b}(\nabla)&=\bigl|T'^{\,\even/\odd}_{g,b}(\nabla)\bigr|
=T'^{\,\even/\odd}_{g,\mu}(\nabla),
\\
R_{\alpha,b}(\nabla)&=\bigl|R_{\alpha,b}(\nabla)\bigr|=R_{\alpha,\mu}(\nabla).
\end{aligned}
\end{equation}

Observe that $\Omega^*(M;E)(0)$ the (generalized) eigen space of
$\Delta^\nabla_{*;g,b}$ corresponding to the eigen value zero is a finite 
dimensional vector space of dimension the multiplicity of $0$. The 
restriction of the symmetric bilinear form induced by $b$ remains 
non-degenerate and defines for each component $\Omega^q(M;E)(0)$ an equivalence 
class of bases. Since $d^\nabla$ commutes with $\Delta^\nabla_{*;g,b}$,  
$\bigl(\Omega^*(M;E)(0),d^\nabla\bigr)$ is a finite dimensional complex. 
When acyclic, i.e.\ $\nabla\in\mathcal F(E)\setminus\Sigma(E)$, denote by
$$
T_\an(\nabla,g,b)(0)\in\mathbb C_*
$$ 
the Milnor torsion associated to the equivalence class of bases induced by $b$.

\subsection {The modified Ray--Singer torsion}\label{SS:5.1}

Let $E\to M$ be a complex vector bundle, 
and let $\e^*\in\Eul^*(M;\mathbb R)$ be a coEuler structure.
Choose a Hermitian structure (fiber metric) $\mu$ on $E$,
a Riemannian metric $g$ on $M$ and $\alpha\in\Omega^{n-1}(M;\mathcal O_M)$
so that $[g,\alpha]=\e^*$, see section~\ref{SS:3.2}.
For $\nabla\in\mathcal F(E)\setminus\Sigma(E)$ consider the quantity
$$
T_\an(\nabla,\mu,g,\alpha):=\bigl(T^{\even}_{g,\mu}(\nabla)
/T^{\odd}_{g,\mu}(\nabla)\bigr)^{-1/2}
\cdot R_{\alpha,\mu}(\nabla)\in\mathbb R_{>0}
$$ 
referred to as the \emph{modified Ray--Singer torsion.} 
The following proposition is a reformulation of one of the main 
theorems in \cite{BZ92}, cf.\ also \cite{BFK01} and \cite{BH04}.

\begin{proposition}\label{P:3}
If $\nabla\in\mathcal F(E)\setminus\Sigma(E)$, then $T_\an(\nabla,\mu,g,\alpha)$
is gauge invariant and independent of $\mu,g,\alpha$.
\end{proposition}

When applied to $\mathbb F_\rho$ the number 
$T_\an^{\e^*}(\rho):=T_\an(\nabla_\rho,\mu,g,\alpha)$
defines a real analytic function
$T_\an^{\e^*}:\Rep^M_0(\Gamma;V)\to\R_{>0}$.
It is natural to ask if $T_\an^{\e^*}$ is the absolute value of a holomorphic function.

The answer is no as one can see on the simplest possible example $M= S^1$ 
equipped with the the canonical coEuler structure $\e^*_0$. 
In this case $\Rep^M(\Gamma;\mathbb C)=\mathbb C\setminus 0$,
and $T_\an^{\e^*_0}(z)=|\frac{(1-z)}{z^{1/2}}|$, cf.\ \cite {BH05}.
However, Theorem~\ref{T:2} in section~\ref{SS:6.1} below provides
the following answer to the question (Q) from the introduction.

\begin{observation}  
If $\e^*$ is an integral coEuler structure, then 
$T_\an^{\e^*}$ is the absolute value of a holomorphic function on 
$\Rep^M_0(\Gamma;V)$ which is the restriction of a rational function on 
$\Rep^M(\Gamma;V)$. For a general coEuler structure $T^{\e^*}_\an$
still locally is the absolute value of a holomorphic function.
\end{observation}

\subsection{Complex Ray--Singer torsion}\label{SS:5.3}

Let $E$ be a complex vector bundle equipped with a non-degenerate symmetric 
bilinear form $b$. Suppose $(g,\alpha)$ is a pair consisting of a Riemannian 
metric $g$ and a differential form $\alpha\in\Omega^{n-1}(M;\mathcal O_M)$ 
with $d\alpha=\ec(g)$. For any $\nabla\in\mathcal F(E)\setminus\Sigma(E)$ 
consider the complex number 
\begin{equation}\label{E:00}
\mathcal{ST}_\an(\nabla,b,g,\alpha)
:=\bigl(T'^{\,\even}_{g,b}(\nabla)/T'^{\,\odd}_{g,b}(\nabla)\bigr)^{-1}
\cdot R_{\alpha,b}(\nabla)^2
\cdot T_\an(\nabla,g,b)(0)^2\in\mathbb C_*
\end{equation}
referred to as the \emph{complex valued Ray--Singer torsion.}\footnote{The idea 
of considering $b$-Laplacians for torsion was brought to the
attention of the  first author by W.~M\"uller \cite{M}.
The second author came to it independently.}

It is possible to provide an alternative definition of 
$\mathcal{ST}_\an(\nabla,b,g,\alpha)$.
Suppose $R>0$ is a positive real number so that the Laplacians $\Delta^\nabla_{q;g,b}$ 
have no eigen values of absolute value $R$.
In this case denote by $\det^R\Delta^\nabla_{q;g,b}$ the regularized product of all eigen
values larger than $R$ w.r.\ to a non-zero Agmon angle disjoint from the spectrum  
$T^{R,\even}_{g,b}$ resp.\ $T^{R,\even}_{g,b}$ the quantities defined by the 
formulae~\eqref{E;16'} with $T^{R,\even/\odd}(\Delta)$ instead of 
$T^{'\,\even/\odd}(\Delta)$.
Consider $\Omega^*(M;E)(R)$ to be the sum of generalized  eigen spaces  
of $\Delta^\nabla_{*;g,b}$  corresponding to eigen values smaller in 
absolute value than $R$. $(\Omega^*(M;E)(R), d^\nabla)$ is a finite 
dimensional complex. As before $b$ remains non-degenerate  and when 
acyclic (and this is the case iff $(\Omega^*(M;E),d^\nabla)$ is acyclic) 
denote by $T_\an(\nabla,g,b)(R)$ the Milnor torsion associated to the equivalence 
class of bases induced by $b$.
It is easy to check that 
\begin{equation}\label{E:000}
\mathcal{ST}_\an(\nabla,b,g,\alpha)
=\bigl(T^{R,\even}_{g,b}(\nabla)/T^{R,\odd}_{g,b}(\nabla)\bigr)^{-1}
\cdot R_{\alpha,b}(\nabla)^2
\cdot T_\an(\nabla,g,b)(R)^2
\end{equation}

\begin{proposition}\label{P:6}
1. $\mathcal{ST}_\an(\nabla,b,g,\alpha)$ is a holomorphic function  
on $\mathcal F(E)\setminus\Sigma(E)$ and the restriction of a meromorphic 
function on $\mathcal F(E)$ with poles and zeros in $\Sigma(E)$.

2. If $b_1$ and $b_2$ are two non-degenerate symmetric bilinear forms which are 
homotopic, then $\mathcal{ST}_\an(\nabla,b_1,g,\alpha)=\mathcal{ST}_\an(\nabla,b_2,g,\alpha)$.

3. If $(g_1,\alpha_1)$ and $(g_2,\alpha_2)$ are two pairs representing the same 
coEuler structure, then 
$\mathcal{ST}_\an(\nabla,b,g_1,\alpha_1)=\mathcal{ST}_\an(\nabla,b,g_2,\alpha_2)$.

4. We have 
$\mathcal{ST}_\an(\gamma\nabla,\gamma b,g,\alpha)=\mathcal{ST}_\an(\nabla,b,g,\alpha)$
for every gauge transformation $\gamma$ of $E$.

5. $\mathcal{ST}_\an(\nabla_1\oplus\nabla_2,b_1\oplus b_2,g,\alpha)= 
\mathcal{ST}_\an(\nabla_1,b_1,g,\alpha)\cdot\mathcal{ST}_\an(\nabla_2,b_2,g,\alpha)$.
\end{proposition}

To check the first part of this proposition, one shows that for $\nabla_0\in\mathcal F(E)$  
one can find $R>0$ and an open neighborhood $U$ of $\nabla_0\in \mathcal F(E)$ such that  
no eigen value of $\Delta^\nabla_{q;g,b}$, $\nabla\in U$, has absolute value $R$. 
The function $\bigl(T^{R,\even}_{g,b}(\nabla)/T^{R,\odd}_{g,b}(\nabla)\bigr)^{-1}$ 
is holomorphic in $\nabla\in U$. Moreover, on $U$ the function $T_\an(\nabla,g,b)(R)^2$ is 
meromorphic in $\nabla$, and holomorphic when restricted to $U\setminus\Sigma(E)$.
The statement thus follows from \eqref{E:000}.

The second and third part of Proposition~\ref{P:6} are derived from formulas for 
$d/dt(\mathcal{ST}_\an(\nabla,b(t),g,\alpha))$ 
resp.\ $d/dt(\mathcal{ST}_\an(\nabla,b,g(t),\alpha)$
which are similar to such formulas for Ray--Singer torsion in the case of a 
Hermitian structure instead of a non-degenerate symmetric bilinear form, cf.\ \cite {BH05}. 
The proof of 4) and 5) require a careful inspection of the definitions.
The full arguments are contained in \cite{BH05}.

As a consequence to each homotopy class of non-degenerate symmetric bilinear forms $[b]$
%(equivalently isomorphisms between $E^*$ and $E$)
and coEuler structure $\e^*$ we can
associate a meromorphic function on $\mathcal F(E)$.  The reader unfamiliar with the basic concepts 
of complex  analytic geometry on Banach/ Frechet manifolds  can  consult \cite{D} and \cite {KM}.   Changing the coEuler structure our 
function changes by multiplication with a non-vanishing holomorphic function as one can 
see from \eqref{E:00}. Changing the homotopy class $[b]$ is actually more subtle. We 
expect however that $\mathcal{ST}$ remains unchanged when the coEuler structure is integral.

Denote by $\Rep^{M,E}(\Gamma;V)$ the union of components of $\Rep^M(\Gamma;V)$
which consists of representations equivalent to holonomy representations of flat 
connections in the bundle $E$. Suppose $E$ admits non-degenerate symmetric bilinear 
forms and let $[b]$ be a homotopy class of such forms. Let $x_0\in M$ be a base point
and denote by $\mathcal G(E)_{x_0,[b]}$ 
the group of gauge transformations which leave fixed $E_{x_0}$ and the class $[b]$.
In view of Proposition~\ref{P:6}, $\mathcal{ST}_\an(\nabla, b, g,\alpha)$ defines
a meromorphic function $\mathcal{ST}_\an^{\e^*,[b]}$ on 
$\pi^{-1}(\Rep^{M,E}(\Gamma;V) \subseteq \mathcal F(E)/\mathcal G_{x_0,[b]}$. Note that
$\pi:\mathcal F(E)/\mathcal G_{x_0,[b]}\to\Rep(\Gamma;V)$
is an principal holomorphic covering of its image which contains $\Rep^{M,E}(\Gamma;V)$.
We expect that the absolute value of this function is the square of 
modified Ray--Singer torsion. The expectation is true when $(E,b)$ satisfies $P_r$ below.

\begin{definition}
The pair $(E,b)$ satisfies \emph{Property $P_r$} if it is the complexification of a 
pair $(E_r,b_r)$ consisting of a real vector bundle $E_r$ and a 
non-degenerate symmetric positive definite $\mathbb R$-bilinear form $b_r$
and the space of flat connections $\mathcal F(E_r)$ 
is a real form of the space $\mathcal F(E)$.
\end{definition}

We summarize this in the following Theorem.

\begin{theorem}\label{T:1}
With the hypotheses above we have.

1.
If $\e^*_1$ and $\e^*_2$ are two coEuler structures then
$$
\mathcal{ST}_\an^{\e^*_1,[b]}=\mathcal{ST}_\an^{\e^*_2,[b]}\cdot 
e^{2([\omega(\nabla,b)],D^{-1}(\e^*_1-\e^*_2))}
$$
with $D:H_1(M;\R)\to H^{n-1}(M;\mathcal O_M)$ the Poincar\'e duality isomorphism.

Suppose that $(E,b)$ satisfies  property $(P_r)$. Then:

2. 
If $\e^* $ is integral then $\mathcal{ST}_\an^{\e^*,[b]}$ is independent of
$[b]$ and descends to a rational function on 
$\Rep^{M,E}(\Gamma;V)$ denoted $\mathcal{ST}_\an^{\e^*}$.

3.
We have
\begin{equation}
\bigl|\mathcal{ST}_\an^{\e^*,[b]}\bigr|=(T_\an^{\e^*}\cdot\pi)^2.
\end{equation}
\end{theorem}

We expect that both 2) and 3) remain true for an arbitrary pair $(E,b)$.

\begin{observation}
Property 5) in Proposition~\ref{P:6} shows that up to multiplication with a root of unity 
the complex Ray--Singer torsion can be defined on all components of
$\Rep^M(\Gamma;V)$, since $F=\oplus_kE$ is trivial for sufficiently large $k$.
\end{observation}

\section{Milnor--Turaev and dynamical torsion}\label{S:6}

\subsection{Milnor--Turaev torsion}\label{SS:6.1}

Consider a smooth triangulation $\tau$ of $M$, and 
choose a collection of orientations $\mathcal O$ of the simplices of
$\tau$. Let $x_0\in M$ be a base point, and set $\Gamma:=\pi_1(M,x_0)$.
Let $V$ be a finite dimensional complex vector space.
For a representation $\rho\in\Rep(\Gamma;V)$, 
consider the chain complex $(C^*_\tau(M;\rho),d^{\mathcal O}_\tau(\rho))$ associated with
the triangulation $\tau$ which computes the cohomology $H^*(M;\rho)$.

Denote the set of simplexes of dimension $q$ by $\mathcal X_q$, and
set $k_i:=\sharp(\mathcal X_i)\cdot\dim(V)$.
Choose a collection of paths $\sigma:=\{\sigma_s\mid s\in\tau\}$ from
$x_0$ to the barycenters of $\tau$ as in section~\ref{SS:3.1}.
Choose an ordering $o$ of the barycenters
and a framing $\epsilon$ of $V$.
Using $\sigma$, $o$ and $\epsilon$ one can identify $C^q_\tau(M;\rho)$ with
$\C^{k_q}$. We obtain in this way a map 
$$
t_{\mathcal O,\sigma,o,\epsilon}:\Rep(\Gamma;V)\to\mathbb D(k_0,\dotsc,k_n)
$$
which sends $\Rep^M_0(\Gamma;V)$ to 
$\mathbb D_\ac(k_0,\dotsc,k_n)$. A look at the explicit definition of
$d^{\mathcal O}_\tau(\rho)$ implies that $t_{\mathcal O,\sigma,o,\epsilon}$ is actually a 
regular map between two algebraic sets. Change of $\mathcal O,\sigma,o,\epsilon$
changes the map $t_{\mathcal O,\sigma,o,\epsilon}$.

Recall that the triangulation $\tau$ determines Euler vector fields $X_\tau$ which 
together with $\sigma$ determine an Euler structure $\e\in\Eul(M;\mathbb Z)$, see
section~\ref{SS:3.1}. Note that the ordering $o$ induces
a cohomology orientation $\mathfrak o$ in $H^*(M;\R)$. In view of the arguments of 
\cite{M66} or \cite{Tu86} one can conclude (cf.\ \cite{BH04}):

\begin{proposition}\label{P:7}
If $\rho\in\Rep^M_0(\Gamma;V)$ different choices of $\tau, \mathcal O, \sigma, o, \epsilon$
provide the same composition $\tau\cdot t_{\mathcal O, \sigma, o, \epsilon}(\rho)$
provided they define the same Euler structure $\e$ and
homology orientation $\mathfrak o$.
\end{proposition}

In view of Proposition~\ref{P:7}
we obtain a well defined complex valued rational function
on $\Rep^M(\Gamma;V)$ called the Milnor--Turaev torsion
and denoted from now on by $\mathcal T_\comb^{\e,\mathfrak o}$.

\begin{theorem}\label{T:2}
1.
The poles and zeros of $\mathcal T_\comb^{\e,\mathfrak o}$ are contained 
in $\Sigma(M)$, the subvariety of representations $\rho$ with $H^*(M;\rho)\neq 0$.

2.
The absolute value of $\mathcal T_\comb^{\e,\mathfrak o}(\rho)$
calculated on $\rho\in\Rep^M_0(\Gamma;V)$ is the modified Ray--Singer 
torsion $T_{\an}^{\e^*}(\rho)$, where $\e^*=\PD(\e)$.

3.
If $\e_1$ and $\e_2$ are two Euler structures then 
$\mathcal T_\comb^{\e_2,\mathfrak o}
=\mathcal T_\comb^{\e_1,\mathfrak o}
\cdot\det_{\e_2-\e_1}$
and $\mathcal T_\comb^{\e,-\mathfrak o}=
(-1)^{\dim V}\cdot\mathcal T_\comb^{\e,\mathfrak o}$ 
where $\det_{\e_2-\e_1}$ is the regular function on $\Rep^M(\Gamma;V)$ 
defined in section~\ref{SS:4.1}.

4.
When restricted to $\Rep^{M,E}(\Gamma;V)$, $E$ a complex vector bundle equipped with a
non-degenerate symmetric bilinear form  $b$ so that $(E,b)$ satisfies Property $P_r$, 
$(\mathcal T_\comb^{\e,\mathfrak o})^2=\mathcal{ST}_\an^{\e^*}$,
where $\e^*=\PD(\e)$.
\end{theorem}

We expect that 4) remains true without any hypothesis.
Parts 1) and 3) follow from the definition and the general properties of
$\tau$, part 2) can be derived from the work of Bismut--Zhang \cite{BZ92} 
cf.\ also \cite{BFK01}, and part 4) is discussed in \cite{BH05}, Remark 5.11.

\subsection{Dynamical torsion}\label{SS:6.2}

Let $X$ be a vector field on $M$ satisfying 
(H), (EG), (L), (MS) and (NCT) from section~\ref{SS:2.3}. 
Choose orientations $\mathcal O$ of the unstable manifolds.
Let $x_0\in M$ be a base point and set $\Gamma:=\pi_1(M,x_0)$.
Let $V$ be a finite dimensional complex vector space.
For a representation $\rho\in\Rep(\Gamma;V)$ 
consider the associated flat bundle $(F_\rho,\nabla_\rho)$, and set
$C^q_X(M;\rho):=\Gamma(F_\rho|_{\mathcal X_q})$, where $\mathcal X_q$
denotes the set of zeros of index $q$. Recall that for 
$x\in\mathcal X$, $y\in\mathcal X$ and every 
homotopy class $\hat\alpha\in\mathcal P_{x,y}$ parallel transport
provides an isomorphism $(\pt^\rho_{\hat\alpha})^{-1}:(F_\rho)_y\to(F_\rho)_x$.
For $x\in\mathcal X_q$ and $y\in\mathcal X_{q-1}$
consider the expression:
%\begin{equation}\label{E:18}
%\delta_X^{\mathcal O}(\rho)_{x,y}:=\sum_{\theta\in\mathcal T(x,y)}
%\epsilon(\theta)(\pt^\rho_\theta)^{-1}
%\end{equation}
\begin{equation}\label{E:18}
\delta_X^{\mathcal O}(\rho)_{x,y}:=\sum_{\hat\alpha\in\mathcal P_{x,y}}
\mathbb I_{x,y}^{X,\mathcal O}(\hat\alpha)(\pt^\rho_{\hat\alpha})^{-1}.
\end{equation}
If the right hand side of \eqref{E:18} is absolutely convergent for all $x$
and $y$ they provide a linear mapping $\delta^{\mathcal O}_X(\rho):C^{q-1}_X(M;\rho)\to
C^q_X(M;\rho)$ which, in view of Proposition~\ref{P:1}(3), makes 
$\bigl(C^*_X(M;\rho),\delta_X^{\mathcal O}(\rho)\bigr)$ a cochain complex.
There is an integration homomorphism
$\Int_X^{\mathcal O}(\rho):\bigl(\Omega^*(M;F_\rho),d^{\nabla_\rho}\bigr)\to
\bigl(C^*_X(M;\rho),\delta_X^{\mathcal O}(\rho)\bigr)$
which does not always induce an isomorphism in cohomology.

Recall that for every $\rho\in\Rep(\Gamma;V)$ 
the composition $\tr\cdot\rho^{-1}:\Gamma\to\C$ factors through conjugacy classes
to a function $\tr\cdot\rho^{-1}:[S^1,M]\to\C$.
Let us also consider the expression
%\begin{equation}\label{E:19}
%P_X(\rho):=\sum_{\hat\theta}(\epsilon(\hat\theta)/p(\hat\theta))
%\tr(\rho(\hat\theta)^{-1})
%\end{equation}
\begin{equation}\label{E:19}
P_X(\rho):=\sum_{\gamma\in[S^1,M]}\mathbb Z_X(\gamma)(\tr\cdot\rho^{-1})(\gamma).
\end{equation}
Again, the right hand side of \eqref{E:19} will in general not converge.

\begin{proposition}\label{P:8}

There exists an open set $U$ in $\Rep^M(\Gamma;V)$, intersecting every irreducible 
component, s.t.\ for any representation $\rho\in U$ we have:

a) The differentials $\delta^{\mathcal O}_X(\rho)$ converge absolutely.

b) The integration $\Int_X^{\mathcal O}(\rho)$ converges absolutely.

c) The integration $\Int_X^{\mathcal O}(\rho)$ induces an isomorphism in cohomology.
 
d) If in addition $\dim V=1$, then 
\begin{equation}\label{E:19'}
\sum_{\sigma\in H_1(M;\Z)/\Tor(H_1(M;\Z))}\,\biggl|\sum_{[\gamma]\in\sigma}
\mathbb Z_X(\gamma)(\tr\cdot\rho^{-1})(\gamma)\biggr|
\end{equation}
converges, cf.~\eqref{E:19}. Here the inner (finite) sum is over all
$\gamma\in[S^1,M]$ which give rise to $\sigma\in H_1(M;\Z)/\Tor(H_1(M;\Z))$.
\end{proposition}

%%%%%%%%%%%%%%%%%
%a) The set $U$ of representations $\rho$ in $\Rep^M_0(\Gamma;V)$ for which the 
%right side of the formula \eqref{E:18} is absolutely convergent 
%has an open interior in any connected component of $\Rep^M(\Gamma;V)$.
%The subset $\Sigma\subset U$ consisting of the representations $\rho\in U$  
%where $\Int_X^{\mathcal O}(\rho)$ does not induce an isomorphism in cohomology 
%is a closed proper complex analytic subset of $U$. 
%b) If $\dim V=1$ the set $U$ of representations $\rho$ in $\Rep^M_0(\Gamma;V)$ for which the 
%right side of the formula 
% \eqref{E:19} 
%is absolutely convergent 
%has an open interior in any connected component of $\Rep^M(\Gamma;V)$. I
%The subset $\Sigma\subset U$ consisting of the representations $\rho\in U$  
%where $\Int_X^{\mathcal O}(\rho)$ does not induce an isomorphism in cohomology 
%is a closed proper complex analytic subset of $U.$ 
%\end{proposition}
%%%%%%%%%%%%%%%%%%%%%%%%%%%%%%%%
This  Proposition is a consequence of exponential growth property (EG) and requires 
(for d)) Hutchings--Lee or Pajitnov results. 
A proof in the case $\dim V=1$ is presented in \cite{BH03'}.
The convergence of \eqref{E:19'} is derived 
from the interpretation of this sum as the Laplace transform of a Dirichlet 
series with a positive abscissa of convergence.

We expect d) to remain true for $V$ of arbitrary dimension.\footnote{Even more, 
we conjecture that \eqref{E:19} converges absolutely on
an open set $U$ as in Proposition~\ref{P:8}.} In this case we make \eqref{E:19}
precise by setting
\begin{equation}\label{E:19''}
P_X(\rho):=\sum_{\sigma\in H_1(M;\Z)/\Tor(H_1(M;\Z))}\,\sum_{[\gamma]\in\sigma}
\mathbb Z_X(\gamma)(\tr\cdot\rho^{-1})(\gamma).
\end{equation}

\begin{observation}
A Lyapunov closed one form $\omega$ for $X$ permits  to consider the family
of regular functions $P_{X;R}$, $R\in\R$, on the variety $\Rep(\Gamma;V)$ 
defined by: 
$$
P_{X;R}(\rho):=\sum_{\hat\theta,-\omega(\hat\theta)\leq R}
(\epsilon(\hat\theta)/p(\hat\theta))\tr(\rho(\hat\theta)^{-1}).
$$ 
If \eqref{E:19'} converges then   
$\lim_{R\to \infty}P_{X;R}$ exists for $\rho$ in an open set of 
representations. We expect that by analytic continuation this can be defined for 
all representations except ones in a proper algebraic subvariety. This is the case 
when $\dim V=1$ or, for $V$ of arbitrary dimension,  when  the  vector field $X$ 
has only finitely many simple closed trajectories.
In this case $\lim_{R\to \infty}P_{X;R}$ has an analytic continuation to a  
rational function on $\Rep(\Gamma;V)$, see section~\ref{S:8} below. 
\end{observation}

As in section~\ref{SS:6.1}, we choose a collection of paths
$\sigma:=\{\sigma_x\mid x\in\mathcal X\}$ from $x_0$ to the zeros of $X$,
an ordering $o'$ of $\mathcal X$, and a framing $\epsilon$ of $V$.
Using $\sigma,o,\epsilon$ we can identify $C^q_X(M;\rho)$ with $\mathbb
C^{k_q}$, where $k_q:=\sharp(\mathcal X_q)\cdot\dim(V)$.
As in the previous section we obtain in this way a holomorphic map
$$
t_{\mathcal O,\sigma,o',\epsilon}:
U\to\hat{\mathbb D}_\ac(k_0,\dotsc,k_n).
$$

An ordering $o'$ of $\mathcal X$ is given by orderings $o'_q$ of 
$\mathcal X_q$, $q=0,1,\dotsc,n$. Two orderings $o'_1$ and $o'_2$ 
are equivalent if $o'_{1,q}$ is obtained from $o'_{2,q}$ by a permutation 
$\pi_q$ so that $\prod_q\sgn(\pi_q)=1$. 
We call an equivalence class of such orderings a \emph{rest point orientation.}
Let us write $\mathfrak o'$ for the rest point orientation determined by $o'$.
Moreover, let $\e$ denote the Euler structure represented by $X$ and $\sigma$, see
Observation~\ref{O:1}.
As in the previous section, the composition $\tau\cdot t_{\mathcal
O,\sigma,o',\epsilon}:U\setminus\Sigma\to\mathbb C_*$ is a holomorphic map
which only depends on $\e$ and $\mathfrak o'$, and will be denoted by
$\tau_X^{\e,\mathfrak o'}$. Consider the holomorphic map 
$P_X:U\to\mathbb C$ defined by formula~\eqref{E:19}. 
The \emph{dynamical torsion} 
is the partially defined holomorphic function
$$
\mathcal T^{\e,\mathfrak o'}_X:=\tau^{\e,\mathfrak o'}_X\cdot e^{P_X}:
U\setminus\Sigma\to\mathbb C_*.
$$

%%%%%%%%%%%%%%%%%%%%%%%%
%%\begin{observation}
%%The absolute convergence of  (\ref {E:19}) is derived in case $\dim V=1$ from the interpretation 
%%of this sum as the Laplace transform of a Dirichlet series with a positive abscissa of convergence.
%A Lyapunov closed one form $\omega$ for $X$ permits to consider the family
%of regular functions on the  variety  on $\R(\Gamma; V),$  $P_{X,R}$  
%$R\in \mathbb R $  defined by :%of regular functions on $\R(\Gamma; V)$,
%%write down   $P_X:\equiv
%%\lim _{R\to \infty} P_{X; R}$ with 
%$$P_{X; R}(\rho):=\sum_{\hat\theta, \omega(\hat\theta)<R}(\epsilon(\hat\theta)/p(\hat\theta))
%\tr(\rho(\hat\theta)^{-1}$$
% We expect that  %as in the case $\dim V=1$ 
%$\lim _{R\to \infty} P_{X; R}$ exists  at least  for $\rho$ in an open set of representations 
%and  by analytic continuation be defined for  all representations except ones in a proper algebraic subvariety.
%This is the case when $\dim V=1$ or,  for $V$ of arbitrary dimension,  when  the  vector field $X$ has only finitely many simple closed trajectories.
%%when
%In this case $\lim _{R\to \infty} P_{X; R}$ has an analytic continuation to a  a rational function 
%on $\R(\Gamma;V), $
%see section 8 below. 
%\end{observation}
%%%%%%%%%%%%%%%%

The following result is based on a theorem of 
Hutchings--Lee and Pajitnov \cite{H02} cf.\  \cite {BH03'}.

\begin{theorem}\label{T:3}  If $\dim V=1$
the partially defined holomorphic function $\mathcal T^{\e,\mathfrak o'}_X$
has an analytic continuation to a rational function equal to 
$\pm\mathcal T_\comb^{\e,\mathfrak o}$.
\end{theorem}

It is hoped  that a  generalization  of Hutchings--Lee and Pajitnov  results 
which will be elaborated in subsequent 
work \cite{BH05b} might led to the proof of the above result for $V$ of arbitrary 
dimension.

\section{Examples}\label{S:7}

\subsection{Milnor--Turaev torsion for mapping tori and twisted Lefschetz zeta function}\label{SS:7.1}

Let $\Gamma_0$ be a group, $\alpha:\Gamma_0\to \Gamma_0$ an isomorphism
and $V$ a complex vector space. Denote by $\Gamma:=\Gamma_0\times_{\alpha}\mathbb Z$ 
the group whose underlying set is $\Gamma_0\times\mathbb Z$ and group operation 
$(g',n) * (g'',m) :=(\alpha^m(g')\cdot g'', n+m)$. 
A representation $\rho:\Gamma\to\GL(V)$ determines a representation 
$\rho_0(\rho):\Gamma_0\to\GL(V)$ the restriction of $\rho$ to $\Gamma_0\times 0$ 
and an isomorphism of $V$, $\theta(\rho)\in\GL(V)$.

Let $(X,x_0)$ be a based  point compact space with $\pi_1(X,x_0)=\Gamma_0$ and 
$f:(X,x_0) \to (X,x_0)$ a homotopy equivalence. For any integer $k$ the map $f$ 
induces the linear isomorphism $f^k: H^k(X; V) \to H^k(X; V)$ and then the standard 
Lefschetz zeta function 
$$
\zeta_f (z):=\frac{\prod_{k\ \even}\det(I-zf^k)}{\prod_{k\ \odd}\det(I-zf^k)}.
$$

More general if $\rho$ is a representation of $\Gamma$ then $f$ and 
$\rho=(\rho_0(\rho), \theta(\rho))$ induce the linear 
isomorphisms $f^k_\rho: H^k(X;\rho_0(\rho)) \to H^k(X; \rho_0(\rho))$ 
and then the $\rho$-twisted  Lefschetz zeta function 
$$
\zeta_f(\rho,z):=\frac{\prod_{k\ \even}\det(I-zf^k_\rho)}{\prod_{k\ \odd}\det(I-zf^k_\rho)}.
$$     
%%%%%%%%%%%%%%%%
%When $\rho$ is the trivial representation the twisted zeta function is the standard zeta function.
%%%%%%%%%%%%%%%%

Let $N$ be a closed connected manifold and $\varphi:N\to N$ 
a diffeomorphism. Without  loss of generality one can suppose 
that $y_0\in N$ is a fixed point of $\varphi$. Define the mapping torus 
$M=N_\varphi$, the manifold obtained from $N\times I$ identifying 
$(x,1)$ with $(\varphi(x),0)$. Let $x_0=(y_0,0)\in M$ be a base point of $M$.
Set $\Gamma_0:=\pi_1(N,n_0)$ and denote by $\alpha:\pi_1(N,y_0)\to \pi_1(N,y_0)$ 
the isomorphism induced by $\varphi$.  We are in the situation considered above 
with $\Gamma=\pi_1(M,x_0)$. The mapping torus structure on $M$ equips $M$ with a 
canonical Euler structure $\e$ and canonical homology orientation $\mathfrak o$. 
The Euler structure $\e$ is defined by any vector field $X$ with $\omega(X)<0$ 
where $\omega:=p^*dt\in\Omega^1(M;\mathbb R)$; all are homotopic. 
The Wang sequence 
\begin{equation}\label{E:Wang}
\cdots\to H^*(M;\mathbb F_\rho)\to H^*(N;i^* (\mathbb F_\rho))
\xrightarrow{\varphi^*_{ \rho} -\id} 
H^*(N;i^*(\mathbb F_\rho) )
\to H^{*+1}(M;\mathbb F_\rho)\to\cdots
\end{equation}
implies $H^*(M;\mathbb F_\rho)=0$ iff $\det(I-\varphi^k_\rho)\neq0$ for all $k$.
The cohomology orientation is derived from the
Wang long exact sequence for the trivial one dimensional real
representation. For details see \cite{BH04}.
We have

\begin{proposition}\label{P:11}
With these notations $\mathcal T_\comb^{\e,\mathfrak o}(\rho)=\zeta_\varphi(\rho,1)$.
\end{proposition}

This result is known cf.\ \cite{BJ96}. A proof can be also derived easily from \cite{BH04}.

\subsection{Vector fields without rest points and Lyapunov cohomology class}\label{SS:7.2}

Let $X$ be a vector field without rest points,
and suppose $X$ satisfies (L) and (NCT).
As in the previous section $X$ defines an Euler structure $\e$.
Consider the expression 
\eqref{E:19}.
By Theorem~\ref{T:3} we have:

\begin{observation}\label{O:5}
With the hypothesis above there exists an open set $U\subseteq\Rep^M(\Gamma;V)$ 
so that \eqref{E:19''} converges,
%$P_X(\rho):=\sum_{\hat\theta}(\epsilon(\hat\theta)/p(\hat\theta))\tr(\rho(\hat\theta)^{-1})$
and $e^{P_X}$ is a well defined holomorphic function on $U$. The function $e^{P_X}$ has 
an analytic continuation to a rational function on $\Rep^M(\Gamma;V)$ equal to 
$\pm\mathcal T^{\e,\mathfrak o}_\comb$. The set $U$ intersects 
non-trivially each connected component of $\Rep^M(\Gamma;V)$.
\end{observation}

\subsection{The Alexander polynomial}\label{SS:7.3}

If $M$ is obtained by surgery on a framed knot, and $\dim V=1$, then
%$\pi_1(M)/[\pi_1(M),\pi_1(M)]=\Z$, 
$\Rep(\Gamma;V)=\C\setminus0$,
and the function
$(z-1)^2  \mathcal T_\comb^{\e,\mathfrak o}$ equals the Alexander polynomial of the knot,
see \cite{Tu02}. Any twisted Alexander polynomial
of the knot can be also recovered from $\mathcal T_\comb^{\e,\mathfrak o}$ for $V$ of 
higher dimension.  One expects  that passing to higher dimensional
representations $\mathcal T_\comb^{\e,\mathfrak o}$ captures even more subtle
knot invariants.

\section{Applications}\label{S:8}

\subsection{The invariant $A^{\e^*}(\rho_1,\rho_2)$}\label{SS:8.1}

Let $M$ be a $V$-acyclic manifold and $\e^*$ a coEuler structure.
Using the modified Ray--Singer torsion we define a $\R/\pi\Z$ valued
invariant (which resembles the Atiyah--Patodi--Singer
spectral flow) for two representations $\rho_1$, $\rho_2$ in the same 
component of $\Rep^M_0(\Gamma;V)$.

By a holomorphic path in $\Rep^M_0(\Gamma;V)$ we understand a holomorphic map
$\tilde\rho:U\to\Rep^M_0(\Gamma;V)$ where $U$ is an open neighborhood 
of the segment of real numbers $[1,2]\times\{0\}\subset\C$ in the 
complex plane. For a coEuler structure $\e^*$
and a holomorphic path $\tilde\rho$ in $\Rep_0^M(\Gamma;V)$ define
\begin{equation}\label{E:01}
\arg^{\e^*}(\tilde\rho)
:=\Re\biggl(2/i\int_1^2\frac{{\partial}(T_{\an}^{\e^*}\circ\tilde\rho)}
{T_{\an}^{\e^*}\circ\tilde\rho}\biggr) \mod\pi.
\end{equation}
Here, for a smooth function $\varphi$ of complex variable $z$, 
$\partial\varphi$ denotes the complex valued $1$-form $(\partial\varphi/\partial z)dz$
and the integration is along the path $[1,2] \times0\subset U$.
Note that

\begin{observation}\label {O:4}
1. 
Suppose $E$ is a complex vector bundle with a non-degenerate bilinear form
$b$, and suppose $\tilde\rho$
is a holomorphic path in $\Rep^{M,E}_0(\Gamma;V)$. Then
$$
\arg^{\e^*}(\tilde\rho)=\arg\Bigl(\mathcal{ST}_\an^{\e^*,[b]}(\tilde\rho(2))
\big/\mathcal{ST}_\an^{\e^*,[b]}(\tilde\rho(1))\Bigr) \mod\pi.
$$
As consequence

2. 
If $\tilde\rho'$ and $\tilde\rho''$ are two holomorphic paths in
$\Rep^M_0(\Gamma;V)$ with $\tilde\rho'(1)=\tilde\rho''(1)$
and $\tilde\rho'(2)=\tilde\rho''(2)$ then
$$
\arg^{\e^*}(\tilde\rho')=\arg^{\e^*}(\tilde\rho'') \mod\pi.
$$

3. 
If $\tilde\rho'$, $\tilde\rho''$ and $\tilde\rho'''$ are three holomorphic
paths in $\Rep^M_0(\Gamma;V)$ with $\tilde\rho'(1)=\tilde\rho'''(1)$,
$\tilde\rho'(2)=\tilde\rho''(1)$ and $\tilde\rho''(2)=\tilde\rho'''(2)$ then
$$
\arg^{\e^*}(\tilde\rho''')=\arg^{\e^*}(\tilde\rho')+\arg^{\e^*}(\tilde\rho'')\mod\pi.
$$
\end{observation}

Observation~\ref{O:4} permits to define a $\R/\pi\Z$ valued numerical
invariant $A^{\e^*}(\rho_1,\rho_2)$ associated to a coEuler 
structure $\e^*$ and two
representations $\rho_1,\rho_2$ in the same connected component of
$\Rep^M_0(\Gamma;V)$.
If there exists a holomorphic path with $\tilde\rho(1)=\rho_1$
and $\tilde\rho(2)=\rho_2$ we set
$$
A^{\e^*}(\rho_1,\rho_2):=\arg^{\e^*}(\tilde\rho)\mod\pi.
$$
Given any two representations $\rho_1$ and $\rho_2$ in the same component of 
$\Rep^M_0(\Gamma;V)$ one can always find a finite collection of holomorphic paths 
$\tilde\rho_i$, $1\leq i\leq k$, in $\Rep^M_0(\Gamma;V)$ so that
$\tilde\rho_i(2)=\tilde\rho_{i+1}(1)$ for all $1\leq i<k$, and such that
$\tilde\rho_1(1)=\rho_1$ and $\tilde\rho_k(2)=\rho_2$.
Then take
$$
A^{\e^*}(\rho_1,\rho_2):=\sum_{i=1}^k\arg^{\e^*}(\tilde\rho_i)\mod\pi.
$$

In view of Observation~\ref{O:4} the invariant is well defined, and if
$\e^*$ is integral it is actually well defined in $\R/2\pi\Z$.
This invariant was first introduced  when the authors were not fully aware of
``the complex Ray--Singer torsion.''
The formula \eqref{E:01} is a more or less obvious
expression of the phase of a holomorphic function in terms of its absolute value,
the Ray--Singer torsion, as positive real valued function.
By Theorem~\ref{T:2} the invariant can  be computed
with combinatorial topology and by section~\ref{S:7}
quite explicitly in some cases. If the representations $\rho_1, \rho_2$
are unimodular then the
coEuler structure is irrelevant. It is interesting to compare this 
invariant to the Atiyah--Patodi--Singer spectral flow; it is not the same but are related.

\subsection{Novikov conjecture}

Let $X$ be  a smooth  vector field which satisfies (H), (L), (MS), (NCT).  
Suppose $\omega$ is a real valued closed one form so that $\omega(X)_x<0$, $x$ 
not a rest point (Lyapunov form). Define the  functions  $I^X_{x,y}:\mathbb R\to\mathbb Z$ 
and $Z^X:\mathbb R\to \mathbb Q$ by  
\begin{equation}
\begin{aligned}
I^{X,\mathcal O}_{x,y}(R):= &\sum _{\hat \alpha,\  \omega (\hat \alpha)<R} 
\mathbb I^{X,\mathcal O}_{x,y} (\hat \alpha)\\
Z^X(R):= &\sum _{\hat \theta,\  \omega (\hat \theta)<R} \mathbb Z^X(\hat\theta) 
\end{aligned}
\end{equation}

Part (a) of the following conjecture was formulated by Novikov for $X=\grad_g\omega$,
$\omega$ a Morse closed one  form  when this vector field satisfies the above properties.

\begin{conjecture}
a) The function $I^{X,\mathcal O}_{x,y} (R)$ has exponential growth.

b) The function $Z^X(R)$ has exponential growth.
\end{conjecture}

Recall that a function $f:\mathbb R \to \mathbb R$ is said to have exponential growth iff
there exists constants $C_1, C_2$ so that  $| f(x) |  <C_1 e^{C_2}$.

As a straight forward consequence of Proposition~\ref{P:8} we have

%%%%%%%%%%%%%%%%%%%%%%%
%\begin{theorem}
%If $M$ is $V$-acyclic for some $V$ then for any vector field $X$ which 
%satisfies (H), (EG), (L), (MS), (NCT)  and $\omega$ a Lyapunov closed one 
%form the above conjecture is true.
%
%The set of vector fields which satisfy (H), (EG), (L), (MS), (NCT)  
%is $C^0$ dense in the space of vector fields which satisfy (H), (L), (MS), (NCT).
%\end{theorem}
%%%%%%%%%%%%%%%%%%%%

\begin{theorem}
a) Suppose $X$ satisfies (H), (MS), (L) and (EG). Then part a) of the
conjecture above holds.

b) Suppose $M$ is $V$-acyclic for some $V$ with $\dim(V)=1$. Moreover, assume
$X$ satisfies (H), (MS), (L), (NTC) and (EG). Then part b) of the conjecture above holds.
\end{theorem}

This result is proved in \cite {BH05}; The $V$-acyclicity in part b) is not
necessary if (EG) is replaced by an apparently stronger assumption (SEG).
Prior to our work Pajitnov has  considered for vector fields which satisfy  (H), 
(L), (MS), (NCT) an additional property, condition  $(\mathfrak C \mathcal Y)$, 
and has verified  part (a) of this conjecture. He has also shown that the vector 
fields which satisfy  (H), (L), (MS), (NCT)  and $(\mathfrak C \mathcal Y)$
are actually $C^0$ dense in the space of vector fields which satisfy (H), 
(L), (MS), (NCT). It is shown in 
\cite{BH05} that Pajitnov vector fields satisfy (EG), and in fact (SEG).

\subsection{A question in dynamics}\label{SS:7.4}

Let $\Gamma$ be a finitely presented group, $V$ a complex vector space and  
$\Rep(\Gamma;V)$ the variety of complex representations.
Consider triples $\underline a:=\{a,\epsilon_-,\epsilon_+\}$  
where $a$ is a conjugacy class of $\Gamma$ and $\epsilon_\pm\in\{\pm1\}$.
Define the rational function $\llet_{\underline a}:\Rep(\Gamma;V)\to\mathbb C$
by
$$
\llet_{\underline a}(\rho)
:=\Bigl(\det\bigl(\id-(-1)^{\epsilon_-}\rho(\alpha)^{-1}\bigr)
\Bigr)^{(-1)^{\epsilon_-+\epsilon_+}}
$$
where $\alpha \in \Gamma$ is a representative of $a.$

Let $(M,x_0)$ be a $V$-acyclic manifold and $\Gamma=\pi_1(M,x_0)$.
Note that $[S^1,M]$ identifies with the conjugacy classes of $\Gamma$.
Suppose $X$ is a vector field satisfying (L) and (NCT).
Every closed trajectory $\hat\theta$ 
gives rise to a conjugacy class $[\hat\theta]\in[S^1,M]$
and two signs $\epsilon_\pm(\hat\theta)$. These signs are obtained
from the differential of the return map in normal direction;
$\epsilon_-(\hat\theta)$ is the parity of the number of real eigenvalues 
larger than $+1$ and $\epsilon_-(\hat\theta)$ is the parity of
the number of real eigenvalues smaller than $-1$.
For a simple closed trajectory, i.e.\ of period $p(\hat\theta)=1$, let
us consider the triple $\underline{\hat\theta}
:=\bigl([\hat\theta],\epsilon_-(\hat\theta),\epsilon_+(\hat\theta)\bigr)$.
This gives a (at most countable) set of triples as in the previous
paragraph.

Let $\xi\in H^1(M;\mathbb R)$ be a Lyapunov cohomology class for $X$.
Recall that for every $R$ there are only finitely many closed
trajectories $\hat\theta$ with $-\xi([\hat\theta])\leq R$.
Hence, we get a rational function $\zeta_R^{X,\xi}:\Rep(\Gamma;V)\to\mathbb C$
$$
\zeta_R^{X,\xi}:=\prod_{-\xi([\hat\theta])\leq R}\llet_{\underline{\hat\theta}}
$$
where the product is over all triples $\underline{\hat\theta}$
associated to simple closed trajectories with $-\xi([\hat\theta])\leq R$.
It is easy to check that formally we have 
$$
\lim_{R\to\infty}\zeta_R^{X,\xi}=e^{P_X}.
$$
It would be interesting to understand in what sense (if any) this can be made precise.
We conjecture that there exists an open set with non-empty interior in each 
component of $\Rep(\Gamma;V)$ on which we have true convergence. 
In fact there exist vector fields $X$ where the sets of triples are 
finite in which case the conjecture is obviously true.  
%%%%
%The case of vector fields with no rest points when $\zeta_{R,\infty}$ is of particular interest.
%%%

\end{document}